\documentclass[11pt]{amsart}
\usepackage{amssymb,amsmath,txfonts,mathrsfs,mathtools,titletoc, tikz, stmaryrd}
\usepackage[alphabetic]{amsrefs}
\usepackage{color}
\usepackage{graphicx}
\usepackage{float}
\usepackage{appendix}

\newtheorem{theorem}{Theorem}[section]
\newtheorem{prop}[theorem]{Proposition}
\newtheorem{lemma}[theorem]{Lemma}
\newtheorem{remark}[theorem]{Remark}

\newtheorem{definition}[theorem]{Definition}
\newtheorem{cor}[theorem]{Corollary}

\newtheorem{conj}[theorem]{Conjecture}

\numberwithin{equation}{section}

\def\pf{{\it Proof:}~}

\begin{document}

\title[The first Dirichlet eigenvalue and the width]{The first Dirichlet eigenvalue and the width}
\author{Guoyi Xu}
\address{Guoyi Xu\\ Department of Mathematical Sciences\\Tsinghua University, Beijing\\P. R. China}
\email{guoyixu@tsinghua.edu.cn}
\date{\today}

\begin{abstract}
For a geodesic ball with non-negative Ricci curvature and mean convex boundary, it is known that the first Dirichlet eigenvalue of this geodesic ball has a sharp lower bound in term of its radius. We show a quantitative explicit inequality, which bounds the width of geodesic ball in terms of the spectral gap between the first Dirichlet eigenvalue and the corresponding sharp lower bound. 
\\[3mm]
Mathematics Subject Classification: 53B20,58C40, 58J05.
\end{abstract}
\thanks{The author was partially supported by NSFC 12141103.}

\maketitle

\titlecontents{section}[0em]{}{\hspace{.5em}}{}{\titlerule*[1pc]{.}\contentspage}
\titlecontents{subsection}[1.5em]{}{\hspace{.5em}}{}{\titlerule*[1pc]{.}\contentspage}
\tableofcontents
\section{Introduction}

If a complete Riemannian manifold $(M^n,  g)$ with $Rc\geq 0$ contains a geodesic line,  then from Cheeger-Gromoll splitting theorem,  we know that $(M^n, g)$ is isometric to $N^{n- 1}\times \mathbb{R}$.  

The above splitting theorem has a local stable version established by Cheeger-Colding \cite{CC-Ann}.  The local splitting theorem says, if there is a line segment $\gamma\subseteq (M^n,  g)$,  then in the mid-point $q$ of $\gamma$,  there exists a small neighborhood $B_q(R)$ which is Gromov-Hausdorff close to $B(R)\subseteq N^{n- 1}\times [-R, R]$; here $N^{n- 1}$ is a length space and $B(R)$ is a metric ball in $N^{n- 1}\times [-R, R]$ with the product metric.  Also see \cite{Xu} for some quantitative discussion.  

The above stable result forms the basis of Cheeger-Colding's theory of Ricci limit spaces.  We are interested in the diameter of $N^{n- 1}$. Note the key technical tool to prove those splitting results is linear growth harmonic function. Another key fact is that the `model'  space is a cylinder. 

We recall the definition of the width. For Riemannian manifold $(M^n, g)$ where $g$ is the Riemannian metric, the width of $M^n$ is defined as follows:
\begin{align}
\mathcal{W}(M^n)= \inf_{f\in \mathrm{Lip}(M^n)}\max_{t\in \mathbb{R}}\mathrm{diam}(f^{-1}(t)), \nonumber 
\end{align}
where $\mathrm{Lip}(M^n)$ is the set of all Lipschitz functions defined on $M^n$, and $\mathrm{diam}(f^{-1}(t))$ is defined by
\begin{align}
\mathrm{diam}(f^{-1}(t))= \max_{x, y\in f^{-1}(t)}d_g(x, y);\nonumber 
\end{align}
and $d_g(x, y)$ the distance between $x, y$ with respect to Riemannian metric $g$.

\begin{figure}[H]
\begin{center}
\includegraphics{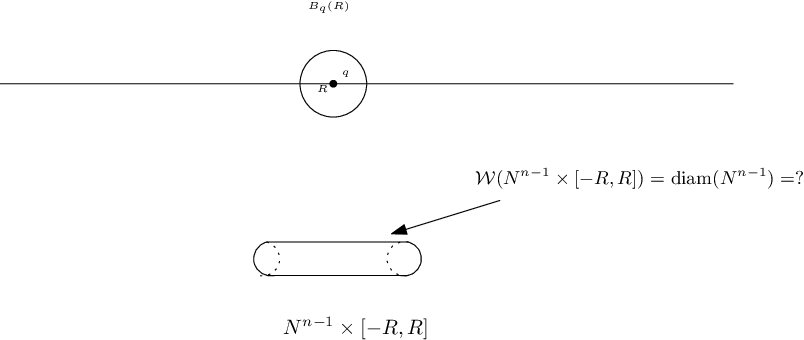}
\caption{Local splitting and `model' cylinder}
\label{figure: splitting}
\end{center}
\end{figure}

In the above context, the width of a geodesic ball depends on two terms. The first term is the Gromov-Hausdorff distance between the geodesic ball and the `model' cylinder, and the second term is the width of the `model' cylinder $N^{n- 1}\times [-R, R]$. Note the width of the `model' cylinder $N^{n- 1}\times [-R, R]$ is not bigger than the diameter of $N^{n- 1}$, when $R$ has much larger scale comparing the scale of $N^{n- 1}$ (see Figure \ref{figure: splitting}). 

We are interested in obtaining an `analytic' characterization of the different scales in a geodesic ball with $Rc\geq 0$. One conjecture of Hang and Wang \cite{HW} as follows partially motivates our study along this direction:
\begin{conj}\label{conj Hang-Wang}
{Assume $\Omega\subseteq \mathbb{R}^2$ is a bounded convex open subset with smooth boundary, then there is a universal constant $c> 0$ such that 
\begin{align}
\mu_1(\Omega)\mathrm{diam}(\Omega)^2- \pi^2\geq c\cdot (\frac{\mathcal{W}(\Omega)}{\mathrm{diam}(\Omega)})^2, \label{width ineq}
\end{align}
where $\mu_1(\Omega)$ is the first non-zero Neumann eigenvalue of $\Omega$. 
}
\end{conj}

Note Payne-Weinberger \cite{PW} proved that $\mu_1(\Omega)\geq \frac{\pi^2}{\mathrm{diam}(\Omega)^2}$ and Hang-Wang \cite{HW} showed the equality holds if and only if $\Omega$ is isometric to an interval of $\mathbb{R}$. The above conjecture points out: the spectral gap between the scaling invariant eigenvalue $\mu_1(\Omega)\mathrm{diam}(\Omega)^2$ and the sharp lower bound $\pi^2$, possibly bounds the square of the scaling invariant width $\frac{\mathcal{W}(\Omega)}{\mathrm{diam}(\Omega)}$. 

In \cite{CMS}, using localization technique developed \cite{CMM} (also see \cite{Klartag}), the inequality between the spectral gap and diameter difference is obtained on manifolds with positive Ricci curvature among more general context.

In this paper, we prove an inequality linking the spectral gap (between the first Dirichlet eigenvalue and the corresponding sharp lower bound) with the width, in similar form as (\ref{width ineq}). 

Recall the inscribed radius of Riemannian manifold $M^n$ is defined as $\displaystyle \mathrm{IR}(M^n)\vcentcolon = \sup_{x\in M^n} d(x,  \partial M^n)$. For the convex region $\Omega\subseteq \mathbb{R}^2$,  Hersch \cite{Hersch} showed $\displaystyle \lambda_1(\Omega)\cdot \mathrm{IR}(\Omega)^2- \frac{\pi^2}{4}\geq 0$. Later, for convex region $\Omega\subseteq \mathbb{R}^n$, Mendez-Hernandez \cite{MH}  (also see Protter \cite{Protter}) gave a proof of the inequality 
\begin{align}
\lambda_1(\Omega)\cdot \mathrm{IR}(\Omega)^2- \frac{\pi^2}{4}\geq \frac{(n- 1)\pi^2}{4}\cdot \Big(\frac{\mathrm{IR}(\Omega)}{\mathrm{diam}(\Omega)}\Big)^2. \label{MH's ineq}
\end{align}

Kasue \cite{Kasue} generalized Hersch's result to manifolds case. More specifically, he proved that for any $n$-dim bounded region $\Omega$ with $Rc\geq 0$ and mean convex boundary, Hersch's inequality also holds. He also showed that the equality case holds if and only if $\Omega$ is isometric to $[0, 2\cdot \mathrm{IR}(\Omega)]\times \tilde{\Omega}$, where $\tilde{\Omega}$ is an $(n-1)$-dim bounded region with $Rc\geq 0$ and mean convex boundary. 

\begin{figure}[H]
\begin{center}
\includegraphics{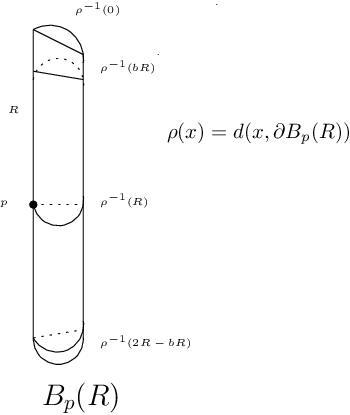}
\caption{Tube type geodesic ball in Theorem \ref{thm width of MCB intr}}
\label{figure: tube}
\end{center}
\end{figure}

In $1999$, among other things, Yang \cite{Yang} used the gradient estimate method (developed by Li-Yau \cite{LY} and Zhong-Yang \cite{ZY}) to present a different proof of Kasue's result. Recently Wang, Zhou and the author \cite{WXZ} prove that $\lambda_1(B_p(R))\cdot R^2- \frac{\pi^2}{16}> 0$ holds for any geodesic ball in non-compact complete Riemannian manifold with $Rc\geq 0$. And the model space in \cite{WXZ} is a tube, see Figure \ref{figure: tube}.  

Motivated by Conjecture \ref{conj Hang-Wang} and the results in \cite{Kasue}, \cite{WXZ}, we try to bound the width of geodesic ball, by the gap between the first Dirichlet eigenvalue and the sharp lower bound. Our main result in this paper is the following theorem. 
\begin{theorem}\label{thm width of MCB intr}
{There is $C(n)> 0$,  such that for any $n$-dim geodesic ball $B_p(R)$ with $Rc\geq 0$ and mean convex boundary (i.e. the mean curvature of $\partial B_p(R)$ with respect to the outforward unit normal vector of $\partial B_p(R)$ is non-negative),  we have  
\begin{align}
\lambda_1(B_p(R))\cdot R^2- \frac{\pi^2}{16}\geq C(n)\Big\{\frac{\mathcal{W}(B_p(R))}{R}\Big\}^{\frac{4(2n+ 1)(8n+ 15)}{n+ 2}}.   \nonumber 
\end{align}
}
\end{theorem}

\begin{remark}\label{rem why we only care about ball}
{Theorem \ref{thm width of MCB intr} partially can be viewed as the stability version of the rigidity result proved in \cite{Kasue}. In this paper we only discuss the case of $\Omega= B_p(R)$. Note the geodesic balls in complete Riemannian manifolds with warped product metric $g= dr^2+ f^2(r)\cdot d\mathbb{S}^{n-1}$ and $Rc(g)\geq 0$, whose center is the origin of the warped product metric (i.e. $r= 0$); have mean convex boundary. 
}
\end{remark}

Recall that by \cite{WXZ}, the inequality $\lambda_1(B_p(R))\cdot R^2- \frac{\pi^2}{16}> 0$ always holds for any geodesic ball in non-compact complete Riemannian manifold with $Rc\geq 0$. And we hope to extend our study to any geodesic balls in non-compact complete Riemannian manifolds in the future. 

The power $\frac{4(2n+ 1)(8n+ 15)}{n+ 2}$ of $\displaystyle \frac{\mathcal{W}(B_p(R))}{R}$ in Theorem \ref{thm width of MCB intr} is not sharp.  Comparing the quantitative isoperimetric inequality (see \cite{FMP}), Faber-Krahn inequality (see \cite{BDV}) and (\ref{MH's ineq}); whether the sharp power is $2$ or not, is an interesting question. 

The whole paper devotes to the proof of Theorem \ref{thm width of MCB intr}. Now we sketch the proof of Theorem \ref{thm width of MCB intr} and the organization of the paper. In Section \ref{sec est of eigen and harmonic}, we firstly observe the key equality (\ref{key ineq}) linking the spectral gap with the integral of the ratio between eigenfunction and model distance function. Then we obtain the estimate of the ratio between eigenfunction and model distance function in the form of the spectral gap (see Proposition \ref{prop ratio is almost 1 for mean convex ball}). Using the $C^0$-estimate of eigenfunction and the integral estimate of model distance function, we obtain the volume comparison for distance function's level sets (see Lemma \ref{lem volume comparison of slice}). 

In Section \ref{sec two est}, following Cheeger-Colding \cite{CC-Ann} and volume comparison of level sets, we get the integral estimate of harmonic function's gradient and Hessian, also $C^0$ estimate of the difference between harmonic function and distance function. One key property about extendability of geodesics in the region is also discussed in quantitative way.

\begin{figure}[H]
\begin{center}
\includegraphics{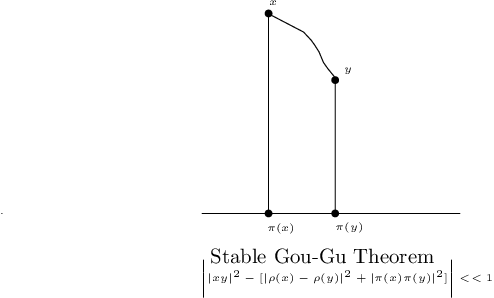}
\caption{Stable-Gou-Gu Theorem}
\label{figure: Stable-Gou-Gu}
\end{center}
\end{figure}

The main result of Section \ref{sec local quantitative metric} is the Stable Gou-Gu Theorem (see Figure \ref{figure: Stable-Gou-Gu} and Theorem \ref{thm Gou-Gu for cylinder type}) for any interior two points of geodesic ball, which are close to each other. Although this technical result (Theorem \ref{thm Gou-Gu for cylinder type}) is possibly well-known to some experts in this field, we provide a detailed argument in this paper; because we can not find the related explicit discussion in the literature.

\begin{figure}[H]
\begin{center}
\includegraphics{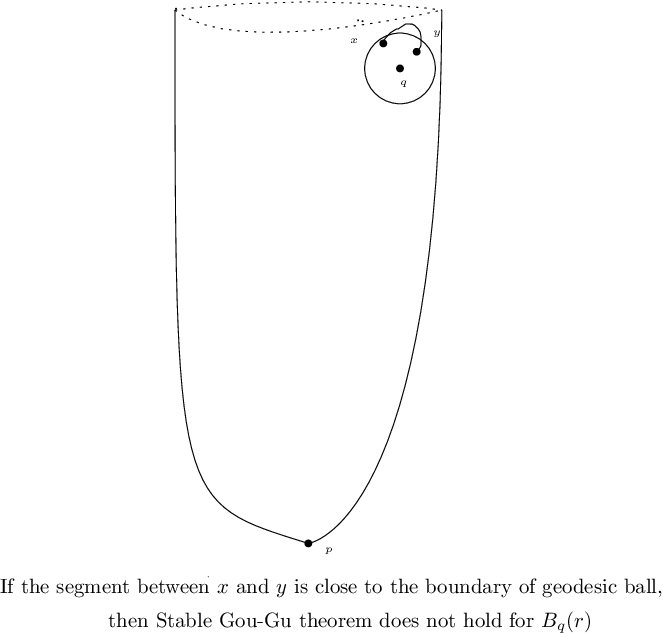}
\caption{Remark}
\label{figure: Ball-near-BD}
\end{center}
\end{figure}

\begin{remark}\label{rem Stable Gou-Gu holds only locally}
{We require that two points are close to each other, because we need that the segment between them does not touch the boundary of geodesic ball. This is a crucial requirement, because we only have global interior integral estimate. When the points or segments are close to the boundary of geodesic ball, the estimate does not hold anymore. We call the distance estimate for such close pair points as `local distance estimate'. 
}
\end{remark}

By modifying the original argument in \cite{CC-Ann} and technical improvement in \cite{CN}, the argument in Section \ref{sec local quantitative metric} and Section \ref{sec width of almost cylinder} provides quantitative argument to obtain the global distance estimate.

Using Theorem \ref{thm Gou-Gu for cylinder type} and the method of getting the global distance estimate from the local distance estimate in \cite{CC-Ann}, in Section \ref{sec width of almost cylinder} we first obtain the estimate for the outer level sets' diameter. Then we use the geodesic ball property, the crucial distance estimates (\ref{assumed Gou-Gu}) and (\ref{geodesic line extendable}), to deduce the width estimate.

Finally we combine the results in Section \ref{sec est of eigen and harmonic}, Section \ref{sec local quantitative metric} and Section \ref{sec width of almost cylinder}, and obtain the width estimate of geodesic ball in term of spectral gap in Section \ref{sec main thm}.

For reader's convenience, we provide a diagram as follows, which is the structure of Theorem \ref{thm width of MCB}(Theorem $1.2$)'s proof.
\begin{figure}[H]
\begin{center}
\includegraphics{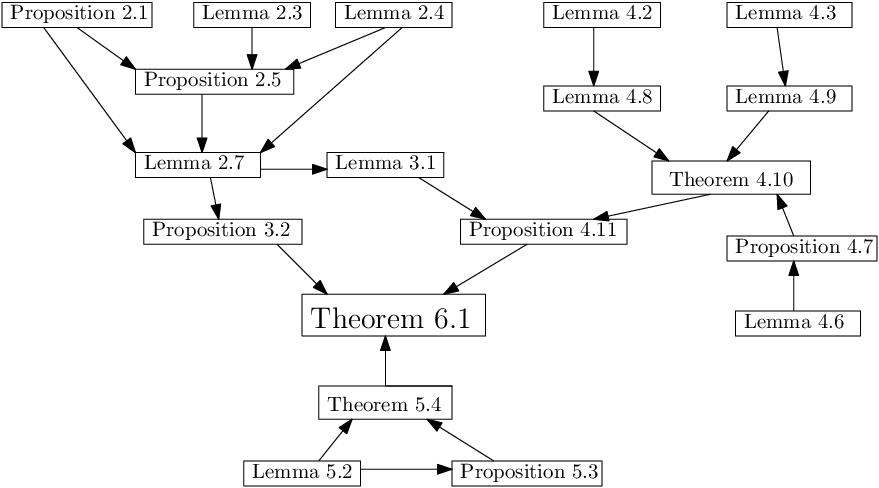}
\caption{The structure of Theorem \ref{thm width of MCB}'s Proof}
\label{figure: Proof}
\end{center}
\end{figure}

\section{Spectral gap bounds the volume difference of level sets}\label{sec est of eigen and harmonic}

In the rest of the paper, unless otherwise mentioned, we assume $u> 0$ is a normalized eigenfunction with respect to the $1$st Dirichlet eigenvalue in $B= B_p(R)$ and $\displaystyle\max_{x\in B} |u(x)|= 1$, where $B$ is a geodesic ball with $Rc\geq 0$ and mean convex boundary. Also we define $\psi(x)= \sin(\frac{\pi}{4R}\rho(x))$ where $\rho(x)= d(x, \partial B_p( R))$. Furthermore, we always assume 
\begin{align}
\epsilon\vcentcolon= \lambda- \frac{\pi^2}{16R^2},\nonumber 
\end{align}
where $\lambda= \lambda_1(B_p(R))$ and $\epsilon\in (0, R^{-2})$. For convenience, for any $C\cdot \epsilon$, where $C> 0$ is some universal constant; we also use $\epsilon$ to denote $C\cdot \epsilon$ in the rest argument. 

In the rest of paper, for any $x, y\in B$, we also use $|xy|$ to denote the distance $d_g(x, y)$ and $\overline{x, y}$ to denote one geodesic segment between $x, y$.

Let $v= \frac{\sin^{-1} u}{\sqrt{\lambda}}$, from \cite{Yang},  we know that 
\begin{align}
|\nabla v|\leq 1.  \label{grad est of v} 
\end{align}

The following key equality of spectral gap is the starting point of our paper.

\begin{prop}\label{prop spectral gap control integ}
{For $u, \psi, \epsilon$ defined as above, we have
\begin{align}
\epsilon&= \frac{\int_B |\nabla (\frac{u}{\psi})|^2\cdot \psi^2}{\int_B u^2}+ \frac{\frac{\pi}{4R}\int_B u^2\cot(\frac{\pi}{4R}\rho)\cdot (-\Delta\rho)}{\int_B u^2}. \label{key ineq}
\end{align}
}
\end{prop}

\begin{remark}\label{rem mean convex implies non-negative terms}
{Define $f(x)= d(p, x)$, then $\rho(x)= R- f(x)$. It is well-known that $\Delta f(x)= H(x)$ in the sense of the distribution, where $H$ is the mean curvature of $\partial B_p(f(x))$ (see \cite[Page $33$]{LiBook}). From $\partial B$ is mean convex and $Rc\geq 0$, we get $\Delta\rho= -\Delta f= -H\leq 0$, where $\Delta \rho$ is again in the sense of the distribution. Therefore, both two terms on the right side of the equality are non-negative. Also the above equation gives another proof of $\lambda_1(\Omega)\cdot \mathrm{IR}(\Omega)^2- \frac{\pi^2}{4}\geq 0$ if we replace $\psi$ by the function $\sin(\frac{\pi}{2\mathrm{IR}(\Omega)}\tilde{\rho}(x))$, where $\tilde{\rho}(x)= d(x, \partial\Omega)$.
}
\end{remark}

\pf
{Define $\eta(x)= \frac{u}{\psi}(x)$, now we have 
\begin{align}
\lambda\int_B u^2&= \int_B |Du|^2= \int_B |\nabla (\psi\cdot \eta)|^2= \int_B \eta^2\cdot |\nabla \psi|^2+ |\nabla \eta|^2\psi^2+ \nabla (\eta^2)\cdot (\psi\nabla \psi) \nonumber \\
&= \int_B |\nabla \eta|^2\psi^2+ \eta^2\cdot |\nabla \psi|^2- \eta^2\nabla (\psi\nabla \psi)= \int_B |\nabla \eta|^2\psi^2- \eta^2\psi\Delta\psi \nonumber \\
&= \frac{\pi^2}{16R^2}\int_B u^2+ \int_B |\nabla (\frac{u}{\psi})|^2\cdot \psi^2+ \frac{\pi}{4R}\int_B u^2\cot(\frac{\pi}{4R}\rho)\cdot (-\Delta\rho). \nonumber 
\end{align}
}
\qed

Firstly we get that the ratio between eigenfunction and distance function is close to $1$ at the center of the ball.
\begin{lemma}\label{lem the value of u at center of ball}
{We have 
\begin{align}
|u(p)- \frac{\sqrt{2}}{2}|\leq R^2\epsilon.  \nonumber 
\end{align}
}
\end{lemma}

\pf
{Assume $\bar{p}\in \partial B_p(R)$, from (\ref{grad est of v}) we get 
\begin{align}
\sin^{-1}u(p)- \sin^{-1}u(\bar{p})\leq \sqrt{\lambda}R\leq  \frac{\pi}{4}+ R^2\epsilon,  \nonumber 
\end{align}
which implies $u(p)\leq \frac{\sqrt{2}}{2}+ R^2\epsilon$. 

On the other hand, assume $u(q)= 1$. Note $d(p, q)\leq R$ we have 
\begin{align}
\sin^{-1}u(q)- \sin^{-1}u(p)\leq \sqrt{\lambda}d(p, q)\leq \frac{\pi}{4}+ R^2\epsilon;  \nonumber 
\end{align}
which implies $u(p)\geq \sin(\frac{\pi}{4}- R^2\epsilon)\geq \frac{\sqrt{2}}{2}- R^2\epsilon$.  
}
\qed

In the rest of the paper, unless otherwise mentioned, the constant $b\in (0, 1)$ is a very small number to be determined later. 

\begin{lemma}\label{lem geodesic ball is rho annulus}
{We have $B_p((1- b)R)= \{bR< \rho < (2- b)R\}$.
}
\end{lemma}

\pf
{Firstly note that $B_p((1- b)R)\subseteq \{bR< \rho < (2- b)R\}$.

Secondly note $\partial B_p(R)$ is a compact set; then for any $x$ with $\rho(x)\in (bR, (2- b)R)$,  there is a point $\hat{x}\in \partial B_p(R)$ with $\rho(x)= d(x, \hat{x})$.  

Note the tangent vector of the geodesic segment $\overline{x, \hat{x}}$ at $\hat{x}$ is perpendicular to the tangent space of $\partial B_p(R)$ at $\hat{x}$, and one writes $\overline{x, \hat{x}}\perp \partial B_p(R)$. This implies 
\begin{align}
\overline{p, \hat{x}}\subseteq \overline{x, \hat{x}},  \quad \quad \text{or}\quad \quad \overline{x, \hat{x}}\subseteq \overline{p, \hat{x}};  \nonumber 
\end{align}
because $\overline{p, \hat{x}}\perp \partial B_p(R)$ at $\hat{x}$ too.  The conclusion follows from the above directly.  
}
\qed

\begin{prop}\label{prop ratio is almost 1 for mean convex ball}
{There is $C_1(n)> 0$ such that if 
\begin{align}
\epsilon R^2\leq C_1(n)b^{n+2}, \label{1st restriction on b}
\end{align}
then there is $C_2(n)> 0$ satisfying
\begin{align}
\sup_{x\in B_p((1- b)R)}|\frac{u}{\psi}(x)- 1|\leq C_2(n)(\epsilon R^2)^{\frac{1}{n+ 2}}b^{-1}. \nonumber 
\end{align}
}
\end{prop}

\begin{remark}\label{rem eigenfunction or harmonic func}
{Because of Proposition \ref{prop ratio is almost 1 for mean convex ball}, in fact we can analyze the metric structure of cylinder type region, using eigenfunctions instead of harmonic functions. However in the context of $Rc\geq 0$, the harmonic function $h$ (see Definition \ref{def of h}) is a canonical approximation to the distance function $\rho$ locally. On the other hand, the eigenfunctions are canonical approximation in the context $Rc\geq (n- 1)$ (see \cite{Colding-sphere-1}, \cite{Colding-sphere-2}, \cite{Petersen}, \cite{Aubry} and \cite{WZ}).  
}
\end{remark}

\pf
{\textbf{Step (1)}. We use the notation $w= \frac{u}{\psi}$ in the rest argument.  From Lemma \ref{lem geodesic ball is rho annulus},  we get $\psi\big|_{B_p((1- b)R)}\geq \sin(\frac{\pi}{4R}b)\geq \frac{1}{2}b$. 

Note $|\nabla v|\leq 1= |\nabla \rho|$, where $|\nabla \rho|= 1$ is in the almost everywhere sense, because $\rho$ is not smooth. We can get 
\begin{align}
u\leq \sin(\sqrt{\lambda}\rho). \nonumber 
\end{align}
Then we have
\begin{align}
\frac{u}{\psi}\leq \frac{\sin(\sqrt{\lambda}\rho)}{\sin(\frac{\pi}{4R}\rho)}\leq 2\sqrt{\lambda}R. \nonumber 
\end{align}

Now we know that
\begin{align}
\sup_{B_p((1- b)R)}|\nabla w|\leq |\frac{|\nabla u|}{\psi}|+ |\frac{u\nabla \psi}{\psi^2}|\leq  \frac{2\sqrt{\lambda}}{b}+ \frac{4\sqrt{\lambda}}{b}\leq C(n)(bR)^{-1}. \label{grad upper bound of w}
\end{align}

\textbf{Step (2)}. From Neumann Poincare inequality for geodesic ball (see \cite{Buser}) and Proposition \ref{prop spectral gap control integ}, note $|u|\leq 1$; we have 
\begin{align}
& \int_{B_p((1- \frac{1}{2}b)R)}|w(x)- \fint_{B_p((1- \frac{1}{2}b)R)} w|^2dx\leq C(n)R^2\int_{B_p((1- \frac{1}{2}b)R)} |\nabla w|^2 \nonumber \\
&\quad\leq C(n)R^2\frac{\int_B |\nabla w|^2 \psi^2}{b^2}\leq  C(n)\epsilon R^2\frac{\int_B u^2}{b^2}\nonumber \\
&\quad\leq C(n)\frac{\epsilon R^2}{b^2}V(B). \label{integ of diff}
\end{align}

Let $A_\eta= \{x\in B_p((1- \frac{1}{2}b)R): |w(x)- \fint_{B_p((1- \frac{1}{2}b)R)} w|\geq \eta\}$, where $\eta= (\epsilon R^2)^{\frac{1}{n+ 2}}b^{-1}$. Then from (\ref{integ of diff}) we obtain
\begin{align}
V(A_\eta)\leq C_3(n)V(B_p(R))\cdot \frac{\epsilon R^2}{b^2\eta^2}. \label{upper bound of volume of bad pts-local}
\end{align}

For any $x\in B_p((1- b)R)$, from Bishop-Gromov Volume Comparison Theorem there is $C_4(n)> 0$ such that if $r_3= C_4(n)\cdot (\frac{\epsilon R^2}{\eta^2 b^2})^{\frac{1}{n}}$, then we have 
\begin{align}
V(B_x(r_3 R))&\geq (\frac{r_3}{2})^n\cdot V(B_x(2R))\geq (\frac{r_3}{2})^n\cdot V(B_p(R))> \frac{C_3(n)V(B_p(R))}{\eta^2}\cdot \frac{\epsilon R^2}{b^2} \nonumber \\
&\geq V(A_\eta). \nonumber 
\end{align}

Now we can find $C_1(n)> 0$, such that if $\epsilon R^2\leq C_1(n)b^{n+ 2}$, then 
\begin{align}
r_3= C_4(n)\cdot (\frac{\epsilon R^2}{\eta^2 b^2})^{\frac{1}{n}}\leq \frac{1}{40}b.\nonumber 
\end{align}
This implies $B_x(r_3 R)\subseteq B_p((1- \frac{1}{2}b)R)$ for any $x\in B_p((1- b)R)$. Therefore for any $x\in B_p((1- b)R)$, there is $\check{x}\in (B_x(r_3R)\backslash A_\eta)\subseteq B_p((1- \frac{1}{2}b)R)$. Now using (\ref{grad upper bound of w}) we get
\begin{align}
&|w(x)- \fint_{B_p((1- \frac{1}{2}b)R)} w|\leq |w(x)- w(\check{x})|+ |w(\check{x})- \fint_{B_p((1- \frac{1}{2}b)R)} w| \nonumber \\
&\leq \sup_{{B_p((1- \frac{1}{2}b)R)}} |\nabla w|\cdot d(x, \check{x})+ \eta\leq \eta+ r_{3}\cdot \frac{C(n)}{b}  \nonumber \\
&\leq \eta+ C_4(n)\cdot (\frac{\epsilon R^2}{\eta^2 b^2})^{\frac{1}{n}}\cdot \frac{1}{b}\leq C(n)(\epsilon R^2)^{\frac{1}{n+ 2}}b^{-1}. \label{fx diff from integ ave}
\end{align}

By Lemma \ref{lem the value of u at center of ball} and (\ref{fx diff from integ ave}), for any $x\in B_p((1- b)R)$ we get
\begin{align}
&|w(x)- 1|\leq |w(x)- \fint_{B_p((1- \frac{1}{2}b)R)} w|+ |w(p)- \fint_{B_p((1- \frac{1}{2}b)R)} w|+ |1- w(p)| \nonumber \\
&\leq C(n)(\epsilon R^2)^{\frac{1}{n+ 2}}b^{-1}+ C(n)R^2\epsilon\leq C(n)(\epsilon R^2)^{\frac{1}{n+ 2}}b^{-1}. \nonumber 
\end{align}
}
\qed

The following lemma provides reason to use the method of Cheeger-Colding in \cite{CC-Ann}.
\begin{lemma}\label{lem volume comparison of slice}
{There is $C(n)> 0$ such that
\begin{align}
\frac{V(\rho^{-1}(bR))- V(\rho^{-1}((2- b)R))}{V(B_p(R))}\leq C(n)\frac{\epsilon R}{b},\nonumber 
\end{align}
where $V(\rho^{-1}(bR)),  V(\rho^{-1}((2- b)R))$ are the $(n- 1)$-dim Hausdorff measure of corresponding sets and $V(B_p(R))$ is the $n$-dim Hausdorff measure of $B_p(R)$.
}
\end{lemma}

\begin{remark}\label{rem CC theory method}
{From the above volume comparison between level sets and Cheeger-Colding \cite[Theorem $4.85$]{CC-Ann},  we in fact get $B_p((1- b)R)$ is close to a cylinder $\rho^{-1}(bR)\times [bR, (2- b)R]$ (which has product metric) in Gromov-Hausdorff distance sense; and we call such set as \textbf{almost cylinder}.  Then we can get the diameter estimate of the level sets in qualitative form. Because we are interested in obtaining the quantitative estimate of the width in form of the spectral gap,   we follow the philosophy of Cheeger-Colding to study the metric structure of $B_p((1- b)R)$ in rest sections. 
}
\end{remark}

\pf
{By Lemma \ref{lem geodesic ball is rho annulus},  let $A= \rho^{-1}(bR, (2- b)R))= B_p((1- b)R)$ in the rest argument.  From Proposition \ref{prop ratio is almost 1 for mean convex ball},  we get 
\begin{align}
\inf_{x\in A}\frac{u}{\psi}\geq \frac{1}{2}.  \label{lower bound of ratio}
\end{align}

From Proposition \ref{prop spectral gap control integ} and (\ref{lower bound of ratio}),  we have
\begin{align}
\epsilon &\geq \frac{\pi}{4R}\frac{\int_B u^2\cot(\frac{\pi}{4R}\rho)(-\Delta\rho)}{\int_B u^2}
\geq \frac{\pi}{8R\int_B u^2}\int_A (\frac{u}{\psi})^2\sin(\frac{\pi}{2R}\rho)(-\Delta\rho) \nonumber \\
&\geq \frac{\pi}{32R\int_B u^2}\int_A \sin(\frac{\pi}{2R}\rho)(-\Delta\rho) .  \label{lower bound by epsilon}
\end{align}

By mean convex property of $\partial B_p(R)$,  we know that $\Delta \rho\leq 0$.  Now by the Divergence Theorem we get
\begin{align}
&\int_A \sin(\frac{\pi}{2R}\rho)(-\Delta\rho) \geq \sin(\frac{b\pi}{2})\cdot \int_A (-\Delta \rho)\geq b\cdot \Big\{V(\rho^{-1}(bR))- V(\rho^{-1}((2- b)R))\Big\}.  \label{upper bound of level set diff}
\end{align}

From (\ref{lower bound by epsilon}) and (\ref{upper bound of level set diff}), we obtain
\begin{align}
\frac{V(\rho^{-1}(bR))- V(\rho^{-1}((2- b)R))}{V(B_p(R))}\leq \frac{\int_A(-\Delta\rho)}{\int_B u^2}\leq C(n)\epsilon Rb^{-1}. \nonumber 
\end{align}
}
\qed

\section{Two type estimates by the volume difference}\label{sec two est}

In this section,  the function $h$ is defined as follows:
\begin{equation}\label{def of h}
\left\{
\begin{array}{rl}
&\Delta h(x)=0 , \quad \quad \quad \quad \quad x\in \rho^{-1}(bR,  (2- b)R),  \\
&h(x)= \rho(x),  \quad \quad \quad \quad \quad x\in \rho^{-1}(bR)\cup \rho^{-1}((2- b)R).  
\end{array} \right.
\end{equation}

\subsection{The integral estimates of harmonic functions}

\begin{lemma}\label{lem grad and Hessian est of h}
{There are $C(n)> 0$ such that 
\begin{align}
&\fint_{\rho^{-1}(bR,  (2- b)R)}|\nabla (h- \rho)|^2\leq C(n)\frac{\epsilon R^2}{b},   \nonumber \\
&\fint_{\rho^{-1}(2bR,  (2- 2b)R)} |\nabla^2 h|^2 \leq \frac{C(n)}{(bR)^2}\cdot \sqrt{\frac{\epsilon R^2}{b}}, \nonumber \\
&\sup_{x\in \rho^{-1}(2bR,  (2- 2b)R)} |h- \rho|(x)\leq C(n)\Big(b\epsilon R^2\Big)^{\frac{1}{n+ 2}}\cdot \frac{R}{b}.  \nonumber 
\end{align}
}
\end{lemma}

\pf
{\textbf{Step (1)}.  Let $A= \rho^{-1}(bR,  (2- b)R)$ and $A_1= \rho^{-1}(2bR,  (2- 2b)R)$.  Firstly by the Maximum principle,  we get $(\rho- h)\big|_A\geq 0$.  Then we have
\begin{align}
\int_A |\nabla (h- \rho)|^2&= \int_A (\rho- h)\Delta (h- \rho)= \int_A (\rho- h)\cdot (-\Delta\rho)\leq \sup_A (\rho- h)\cdot \int_A (-\Delta \rho) \nonumber \\
&\leq (2- b)R\cdot \Big[V(\rho^{-1}(bR))- V\big(\rho^{-1}((2- b)R)\big)\Big].  \nonumber 
\end{align}

Note $B_p((1- b)R)\subseteq A$,  from Lemma \ref{lem volume comparison of slice} we have 
\begin{align}
\fint_A |\nabla (h-\rho)|^2&\leq \frac{\int_A |\nabla (h-\rho)|^2}{V(B_p(R))}\cdot (\frac{R}{(1- b)R})^n \nonumber \\
&\leq C(n)\frac{(2- b)R\cdot \Big[V(\rho^{-1}(bR))- V\big(\rho^{-1}((2- b)R)\big)\Big]}{V(B_p(R))}\leq C(n)\frac{\epsilon R^2}{b}.  \label{grad est of h}
\end{align}

\textbf{Step (2)}.  From Bochner formula and $Rc\geq 0$,  we get $|\nabla ^2 h|^2\leq \frac{1}{2}\Delta (|\nabla h|^2)$.  Therefore choosing cut-off function $\varphi$ with $\varphi\big|_{A_1}= 1$ and $\varphi|_{B- A}= 0$, we obtain
\begin{align}
&\int_{A_1} |\nabla^2 h|^2\leq \int_A |\nabla^2 h|^2\cdot \varphi \leq \frac{1}{2}\int_A\varphi\cdot  \Delta (|\nabla h|^2- 1)\leq \frac{1}{2}\int_A\Delta\varphi\cdot (|\nabla h|^2- 1) \nonumber \\
&\leq \frac{1}{2}\sup_A |\Delta \varphi| \int_A  |\nabla (h- \rho)|\cdot (|\nabla h|+ 1)\leq \frac{C(n)}{(bR)^2}\int_A |\nabla (h- \rho)|\cdot (|\nabla (h- \rho)|+ 2) \nonumber \\
&\leq \frac{C(n)}{(bR)^2}\Big\{\int_A |\nabla (h- \rho)|^2+ \sqrt{V(A)\cdot \int_A |\nabla (h- \rho)|^2}\Big\}.\nonumber 
\end{align}

By volume comparison Theorem and the above,  also using (\ref{grad est of h}),  we have
\begin{align}
\fint_{A_1}|\nabla^2 h|^2\leq \frac{C(n)}{(bR)^2}\Big\{\fint_A |\nabla (h- \rho)|^2+ \sqrt{\fint_A |\nabla (h- \rho)|^2}\Big\}\leq \frac{C(n)}{(bR)^2}\cdot \sqrt{\frac{\epsilon R^2}{b}}.  \nonumber 
\end{align}

\textbf{Step (3)}. Note $A\subseteq B_p(R)$, from the domain monotonicity of Dirichlet eigenvalues (see \cite[Page $17$]{Chavel}); we know that
\begin{align}
\lambda_1(A)\geq \lambda_1(B_p(R))> \frac{\pi^2}{16R^2}. \nonumber 
\end{align}
Combining the result from Step (1),  we have  
\begin{align}
\fint_{A} |h- \rho|^2\leq \lambda_1(A)^{-1}\fint_{A}|\nabla (h- \rho)|^2 \leq C(n)\frac{\epsilon R^2}{b} R^2.  \label{h integral upper bound} 
\end{align}

Let $A_2= \rho^{-1}(\frac{3}{2}bR, (2- \frac{3}{2})bR)$, from Cheng-Yau's gradient estimate for positive harmonic function $h$, we get 
\begin{align}
\sup_{A_2}|\nabla (h- \rho)|\leq \sup_{A_2}(|\nabla h|+ 1)\leq \frac{C(n)}{bR}\sup_{A}h(x)+ 1\leq \frac{C(n)}{b}.  \label{grad est of h pt-wise}
\end{align}

For $\eta> 0$ to be determined later, we define 
\begin{align}
W= h- \rho, \quad \quad \quad \Omega_\eta= \{z\in A_2: |W(z)|\leq \eta\}.  \nonumber 
\end{align}
Let $r_\eta= C(n)R\cdot \Big(\frac{\epsilon R^2}{b}\frac{R^2}{\eta^2}\Big)^{\frac{1}{n}}$, for $x\in A_1$ we have $B_x(r_\eta)\subseteq A_2$. From (\ref{h integral upper bound}) we can get
\begin{align}
\frac{V(B_x(r_\eta))}{V(A)}&\geq (\frac{r_\eta}{2R})^n> C(n)\frac{\epsilon R^2}{b} R^2 \eta^{-2}\geq \eta^{-2}\fint_{A}|h- \rho|^2 \geq \frac{V(A_2- \Omega_\eta)}{V(A)}.  \nonumber 
\end{align}
Therefore $V(B_x(r_\eta))> V(A_2- \Omega_\eta)$.

Note $B_x(r_\eta)\subseteq A_2$, then we get that $B_x( r_\eta)\cap \Omega_\eta\neq \emptyset$. Choose $\check{x}\in B_x( r_\eta)\cap \Omega_\eta$, now by (\ref{grad est of h pt-wise}),  we get
\begin{align}
|W(x)|\leq |W(x)- W(\check{x})|+ |W(\check{x})|\leq \eta+ r_\eta\cdot \sup_{A_2}|\nabla W| \leq \eta+ \frac{C(n)R}{b}\cdot \Big(\frac{\epsilon R^2}{b}\frac{R^2}{\eta^2}\Big)^{\frac{1}{n}}.  \nonumber 
\end{align}

Now choose $\eta= \frac{C(n)R}{b}\cdot \Big(\frac{\epsilon R^2}{b}\frac{R^2}{\eta^2}\Big)^{\frac{1}{n}}$,  which implies $\eta= C(n)\Big(b\epsilon R^2\Big)^{\frac{1}{n+ 2}}\cdot \frac{R}{b}$.  The conclusion follows from the above.
}
\qed

\subsection{The error estimates of extendability}

The next proposition says that the gradient curves of distance function are almost extendable.
\begin{prop}\label{prop grad curv of rho is extendable}
{For $x\in \rho^{-1}(bR)$ and $t\in [b, 2-b]$, there is $z\in \rho^{-1}(tR)$ such that 
\begin{align}
d(x, z)- (t- b)R\leq C(n)R\cdot (\frac{\epsilon R^2}{b})^{\frac{1}{n- 1}}. \nonumber 
\end{align}
}
\end{prop}

\begin{figure}[H]
\begin{center}
\includegraphics{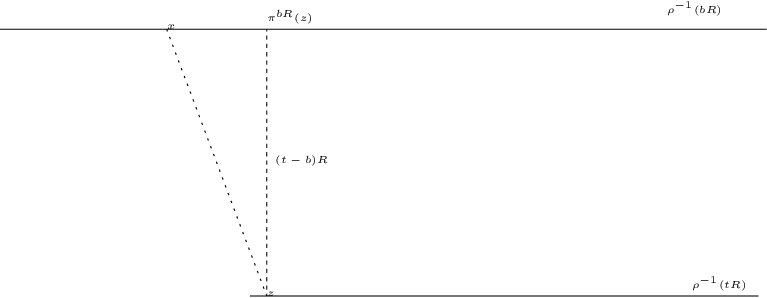}
\caption{gradient curve of $\rho$ is almost extendable}
\label{figure: 2}
\end{center}
\end{figure}

\pf
{\textbf{Step (1)}. We firstly show the conclusion for $t= 2- b$: for any $x\in \rho^{-1}(bR)$  there is $y\in \rho^{-1}((2- b)R)$ such that 
\begin{align}
d(x, y)- (2- 2b)R\leq C(n)R\cdot (\frac{\epsilon R^2}{b})^{\frac{1}{n- 1}}. \label{conclusion for t=2- b}
\end{align}

Let $\zeta= C(n)\cdot (\frac{\epsilon R^2}{b})^{\frac{1}{n- 1}}$. By contradiction, assume the conclusion does not hold for $x_0\in \rho^{-1}(bR)$, which implies
\begin{align}
d(x_0, y)- (2- 2b)R> \zeta R, \quad \quad \quad \quad \forall y\in \rho^{-1}((2- b)R). \nonumber 
\end{align}
Then for any $z\in B_{x_0}(\frac{\zeta R}{2})$ and $y\in \rho^{-1}((2- b)R)$, we have
\begin{align}
d(z, y)- ((2- b)R- \rho(z))&\geq d(x_0, y)- (2- 2b)R- d(z, x)- |bR- \rho(z)| \nonumber \\
&> \zeta R- \frac{\zeta}{2}R- \frac{\zeta}{2}R= 0. \nonumber 
\end{align}

The above inequality implies that $\rho^{-1}((2- b)R)\subseteq \gamma(\rho^{-1}(r R)- B_{x_0}(\frac{\zeta}{2}R))$ for any $r< 2- b$, where $\gamma$ is the gradient curve of $\rho$. 

\textbf{Step (2)}. Note $B_{x_0}(\frac{\zeta}{2}R)\subseteq \rho^{-1}(bR- \frac{\zeta}{2}R, bR+ \frac{\zeta}{2}R)$, by the Co-Area formula we get that there is $r_0\in [b- \frac{\zeta}{2}, b+ \frac{\zeta}{2}]$ such that 
\begin{align}
V(\rho^{-1}(r_0R)\cap B_{x_0}(\frac{\zeta}{2}R))&\geq \frac{1}{\zeta R}V(B_{x_0}(\frac{\zeta}{2}R))\geq \frac{\zeta^{n- 1}}{4^nR}V(B_{x_0}(2R)) \nonumber \\
&\geq \frac{\zeta^{n- 1}}{4^nR}V(B_p(R)), \label{lower bound of the slice}
\end{align}
where the last but one inequality above follows from Bishop-Gromov volume comparison Theorem.

By $\Delta\rho\leq 0$, we know the volume element of the slice $\rho^{-1}(t)$ is decreasing along $t$; because $\Delta\rho$ is the directional derivative of the volume element of the slice $\rho^{-1}(t)$ along $t$ (see \cite[Page $4$ and Page $33$]{LiBook}). Therefore by result in Step (1) and $\frac{b}{2}\leq r_0$, we get
\begin{align}
V(\rho^{-1}((2- b)R))&\leq V(\gamma(\rho^{-1}(r_0R)- B_{x_0}(\frac{\zeta R}{2})))\leq V(\rho^{-1}(r_0R)- B_{x_0}(\frac{\zeta R}{2})) \nonumber \\
&\leq V(\rho^{-1}(r_0R))- V(\rho^{-1}(r_0R)\cap B_{x_0}(\frac{\zeta R}{2})) \nonumber \\
&\leq V(\rho^{-1}(\frac{b R}{2}))- C(n)\zeta^{n- 1}R^{-1}V(B_p(R)). \nonumber 
\end{align}

Therefore by Lemma \ref{lem volume comparison of slice}, we obtain
\begin{align}
\zeta^{n- 1}R^{-1}&\leq C(n)\frac{V(\rho^{-1}(\frac{b}{2}R))- V(\rho^{-1}((2- b)R))}{V(B_p(R))} \nonumber \\
&\leq C(n)\frac{V(\rho^{-1}(\frac{b}{2}R))- V(\rho^{-1}((2- \frac{b}{2})R))}{V(B_p(R))} \nonumber \\
&\leq C(n)\cdot \frac{\epsilon R}{2^{-1}b} . \nonumber 
\end{align}
This contradicts the choice of $\zeta$ at the beginning of the proof. Hence the conclusion (\ref{conclusion for t=2- b}) holds.

Finally we show the conclusion for $t\in [b, 2-b]$. For $x, y$ as in (\ref{conclusion for t=2- b}), we choose $z\in \overline{x, y}\cap \rho^{-1}(tR)$, then
\begin{align}
|zx|&= |xy|- |zy|\leq (2- 2b)R+ C(n)R\cdot (\frac{\epsilon R^2}{b})^{\frac{1}{n- 1}}- |\rho(z)- \rho(y)| \nonumber \\
&= (t- b)R+ C(n)R\cdot (\frac{\epsilon R^2}{b})^{\frac{1}{n- 1}}. \nonumber 
\end{align}
}
\qed

\section{Stable Gou-Gu Theorem with integral estimates}\label{sec local quantitative metric}

In \cite{CC-Ann}, there is an original, concise argument explaining how to obtain the local distance estimate. In \cite{CN}, there is a detailed argument around estimating the distance by the integral of suitable terms, which provides effective ways to get quantitative estimate locally. In \cite{Xu} and \cite{XZ}, following the strategy of \cite{CC-Ann} and using the technical tools developed in \cite{CN}, we essentially proved the following almost Gou-Gu Theorem:
\begin{theorem}\label{thm almost Gou-Gu}
{There are $C(n)> 0$, such that for any $|xy|\leq r$, we can find $\hat{x}, \hat{y}$ satisfying 
\begin{align}
&|x\hat{x}|+ |y\hat{y}|\leq  C(n)r\cdot \zeta(\sup_{B_x( 4r)}\sqrt{\frac{|h- \rho|}{r}},  \fint_A H^2,  \fint_A G^2), \nonumber \\
&\Big|d(x, y)- \sqrt{|\rho(x)- \rho(y)|^2+ |\pi(\hat{x})\pi(\hat{y})|^2}\Big| \nonumber \\
&\leq C(n)r\cdot \zeta(\sup_{B_x( 4r)}\sqrt{\frac{|h- \rho|}{r}},  \fint_A H^2,  \fint_A G^2), \nonumber 
\end{align}
where $\zeta$ is some function satisfying $\displaystyle \lim_{\sum\limits_{i= 1}^3|x_i|\rightarrow 0}\zeta(x_1, x_2, x_3)= 0$.
}
\end{theorem}

\begin{figure}[H]
\begin{center}
\includegraphics{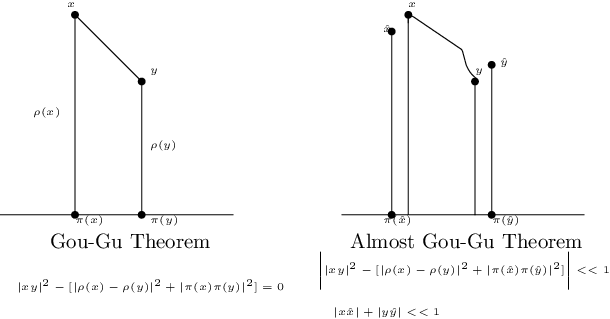}
\caption{Gou-Gu Theorem and Almost Gou-Gu Theorem}
\label{figure: Almost-Gou-Gu}
\end{center}
\end{figure}

Basically, the Almost Gou-Gu Theorem says that the distance between two close points $x, y$ can be controlled by their distance to the base level set $\rho^{-1}(t)$, and the distance between the projection of two  `good' points $\hat{x}, \hat{y}$ on $\rho^{-1}(t)$, where $\hat{x}, \hat{y}$ are close to $x, y$ and the corresponding integral estimates hold for `good' points (see Figure \ref{figure: Almost-Gou-Gu}). 

For $x\in B$, we firstly define $\xi_x$ as the gradient curve of $-\rho$ starting from $x$ and $\gamma_x$ is the gradient curve of $\rho$ starting at $x$. Because $|\nabla \rho|= 1$ almost everywhere, we know $\xi_x, \gamma_x$ is well-defined and unique for almost every $x\in B$; where the measure is the Hausdorff measure determined by Riemannian metric $g$ on $B$. In the rest of the paper, unless otherwise mentioned, all inequalities about distance estimates hold for almost every corresponding point on $B$. 

 In this section, for $0< r<< R$ we assume
\begin{align}
&x\in B,  \quad \quad \rho(x)= b_2,  \quad \quad \frac{1}{2}r\leq b_2- b_0\leq \frac{3}{2}r,  \quad \quad \pi^{b_0}(x)=\xi_x\cap \rho^{-1}(b_0),  \nonumber \\
&|xy|\leq r, \quad \quad \quad \rho(y)- b_0\geq \frac{1}{2}r; \quad \quad \quad b_0\geq 20r, \quad \quad b_2\leq 2R- 20r .\nonumber 
\end{align}

We study the relation between the distance function and the projection map onto $\rho^{-1}(b_0)$. One general philosophy of distance estimate and derivative estimate of distance function on manifolds is, to reduce the estimates to the line integral gradient estimate of error functions between the distance function $\rho$ and harmonic function $h$ and the Hessian estimate of $h$, where $h$ is defined in (\ref{def of h}).

In the rest of this section, for simplicity, we use the notation $G= |\nabla (\rho- h)|, H= |\nabla^2 h|$ (where $G, H$ means gradient and Hessian respectively) and define
\begin{align} 
\mathscr{G}^{b_0}(x)= (\pi^{b_0}(x), \rho(x))\in \rho^{-1}(b_0)\times \mathbb{R};\nonumber 
\end{align}
where $\rho^{-1}(b_0)\times \mathbb{R}$ is endowed with the product distance. 

\begin{remark}\label{rem measurable property}
{Note $\rho$ is a Lipschitz function, from Rademacher's Theorem, we know that $\nabla \rho$ is bounded measurable. From the theory of ODE, we know $\mathscr{G}^{b_0}$ is a measurable map.  Similarly, other maps in later argument are related to $\rho$ or $\nabla \rho$, hence are all measurable maps.
}
\end{remark}

\subsection{Distance estimate by line integral of difference}

We firstly estimate the distance $d(x,y)$ in term of line integrals.

\begin{lemma}\label{lem property 2 of G-H appr}
{Define $\displaystyle l_t:= d\big(\xi_x(t), \tilde{\xi}_y(t)\big)$ for $t\in [0, 2r]$, where
\begin{align*}
\tilde{\xi}_y(t)=\xi_y(\frac{\rho(y)- b_0}{\rho(x)- b_0} t),\ \ \  \tau_t(s)= \overline{\xi_x(t), \tilde{\xi}_y(t)}(s) \text{ for } s\in [0, l_t].
\end{align*}  
Then we have
\begin{align}
&\quad \quad \Big||xy|^2-  |\mathscr{G}^{b_0}(x)\mathscr{G}^{b_0}(y)|^2\Big| \nonumber\\
&\leq 16r\cdot(\int_0^{\rho(y)- b_0} G(\xi_y(t))dt+ \int_0^{\rho(x)- b_0}G(\xi_x(t))dt) \nonumber \\
&\quad  +4r\cdot \sup_{B_x( 4r)}|h- \rho|+ 64r\cdot \int_0^{\rho(x)- b_0} dt\int_{0}^{l_t} H(\tau_t(s)) ds.\nonumber
\end{align}
}
\end{lemma}

\begin{figure}[H]
\begin{center}
\includegraphics{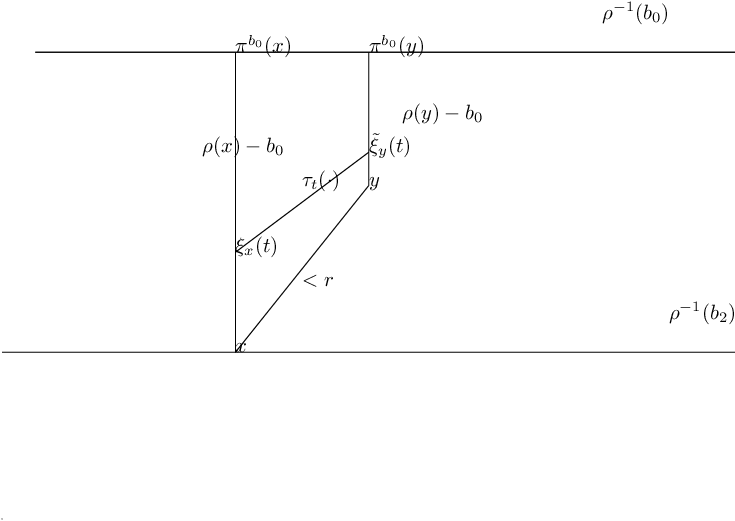}
\caption{Distance estimate of $|xy|$}
\label{figure: lemma2.1}
\end{center}
\end{figure}

\pf
{In the rest of this proof,  for simplicity we use $\pi(x)$ to denote $\pi^{b_0}(x)$. 

Assume $\displaystyle \sup_{x\in B_x( 4r)}|h- \rho|(x)= \delta$. Define $\displaystyle \beta= \frac{\rho(y)- b_0}{\rho(x)- b_0}$ and note $\beta\in [0, 2]$. Since $l_t$ is almost everywhere differentiable on $[0, 2r]$ with respect to $t$, from the extension of the first variation formula (see Theorem \ref{thm general of 1st vari form} in Appendix or \cite{Liu}), we have
\begin{align}
\frac{d}{dt}l_t&= \beta\langle \nabla \rho, \tau_t'\rangle(l_t)- \langle \nabla \rho, \tau_t'\rangle(0) = (I)+ (II)+ \frac{(\beta- 1)^2}{l_t}t, \nonumber \\
&(I)= \beta\langle \nabla (\rho- h), \tau_t'\rangle(l_t)- \langle \nabla (\rho- h), \tau_t'\rangle(0), \nonumber \\
&(II)= \beta[(h\circ \tau_t)'(l_t)- u_t'(l_t)]+ [(h\circ \tau_t)'(0)- u_t'(0)], \nonumber \\
&u_t(s)= \rho(x)- t- \frac{(\beta- 1)t+ \rho(x)- \rho(y)}{l_t}s. \nonumber 
\end{align}

The above inequality yields
\begin{align}
[l_t^2- (\beta- 1)^2t^2]'= 2l_t\cdot (I)+ 2l_t\cdot (II). \nonumber 
\end{align}

Taking integral with respect to $t$ from $0$ to $\rho(x)- b_0$ in the above, note $l_t\leq 8r$, we get
\begin{align}
&\Big|d(x, y)^2- \Big[d(\pi(x),  \pi(y))^2+ |\rho(x)- \rho(y)|^2\Big]\Big| \leq \int_0^{\rho(x)- b_0} |[l_t^2- (\beta- 1)^2t^2]'|dt \nonumber \\
&\leq 16r\int_0^{\rho(x)- b_0} |(I)|dt+ \int_0^{\rho(x)- b_0} 2l_t\cdot |(II)|. \nonumber 
\end{align}

Direct computation yields
\begin{align}
\int_0^{\rho(x)- b_0} |(I)|dt\leq \int_0^{\rho(y)- b_0} |\nabla (\rho- h)|(\xi_y(s))ds+ \int_0^{\rho(x)- b_0} |\nabla (\rho- h)|(\xi_x(s))ds. \nonumber 
\end{align}

By the Mean Value Theorem, there is $\zeta\in (0, l_t)$ such that
\begin{align}
&|(h\circ \tau_t)'(\zeta)- \frac{\beta- 1}{l_t}t|= |\frac{[(h\circ \tau_t)- u_t](l_t)+ [(h\circ \tau_t)- u_t](0)}{l_t}| \nonumber \\
&= |\frac{(h- \rho)(\xi_y(\beta t))- (h- \rho)(\xi_x(t))}{l_t}|\leq \frac{2\delta}{l_t}
\end{align}

Now note $\beta\in (0, \beta_0]$, we get
\begin{align}
|(II)|&= |(\beta- 1)[(h\circ \tau_t)'(l_t)- \frac{\beta- 1}{l_t}t]+ \int_0^{l_t} \frac{d^2}{ds^2}(h\circ \tau_t)(s)ds| \nonumber \\
&\leq |\beta- 1|\cdot |(h\circ \tau_t)'(\zeta)- \frac{\beta- 1}{l_t}t|+ [|\beta- 1|+ 1]\cdot \int_0^{l_t} |\nabla^2 h|(\tau_t(s))ds \nonumber \\
&\leq \frac{2|\beta- 1|\delta}{l_t}+ 4\int_0^{l_t}|\nabla^2 h|(\tau_t(s))ds. \nonumber 
\end{align}

Combining all the above, we obtain
\begin{align}
&\Big|d(x, y)^2- \Big[d(\pi(x),  \pi(y))^2+ |\rho(x)- \rho(y)|^2\Big]\Big| \leq 16r\int_0^{\rho(x)- b_0} |(I)|dt+ \int_0^{\rho(x)- b_0} 2l_t\cdot |(II)| \nonumber \\
&\leq 16r\cdot(\int_0^{\rho(y)- b_0} |\nabla (\rho- h)|(\xi_y(t))dt+ \int_0^{\rho(x)- b_0} |\nabla (\rho- h)|(\xi_x(t))dt) \nonumber \\
&+ 4\delta r+ 64r\int_0^{\rho(x)- b_0} dt\int_{0}^{l_t} |\nabla^2 h|(\tau_t(s)) ds.\nonumber
\end{align}
}
\qed

Assume $b_0- b_1= 10\delta r$ in the rest of the paper.

\begin{lemma}\label{lem Gou-Gu for two balls}
{For $x_1\in B_x(\delta r), x_2\in B_{\pi^{b_0}(x)}(\delta r)$ where $\delta<< 1$, we have 
\begin{align}
&\quad \quad \Big|d(x_1, x_2)^2- d(\mathscr{G}^{b_1}(x_1),  \mathscr{G}^{b_1}(x_2))^2\Big| \nonumber\\
&\leq 16r\cdot(\int_0^{\rho(x_2)- b_1} G(\xi_{x_2}(t))dt+ \int_0^{\rho(x_1)- b_1}G(\xi_{x_1}(t))dt) \nonumber \\
&\quad  +8r\cdot \sup_{B_x( 4r)}|h- \rho|+ 64r\cdot \int_0^{\rho(x_1)- b_1} dt\int_{0}^{l_t} H(\tau_t(s)) ds.\nonumber
\end{align}
}
\end{lemma}

\begin{figure}[H]
\begin{center}
\includegraphics{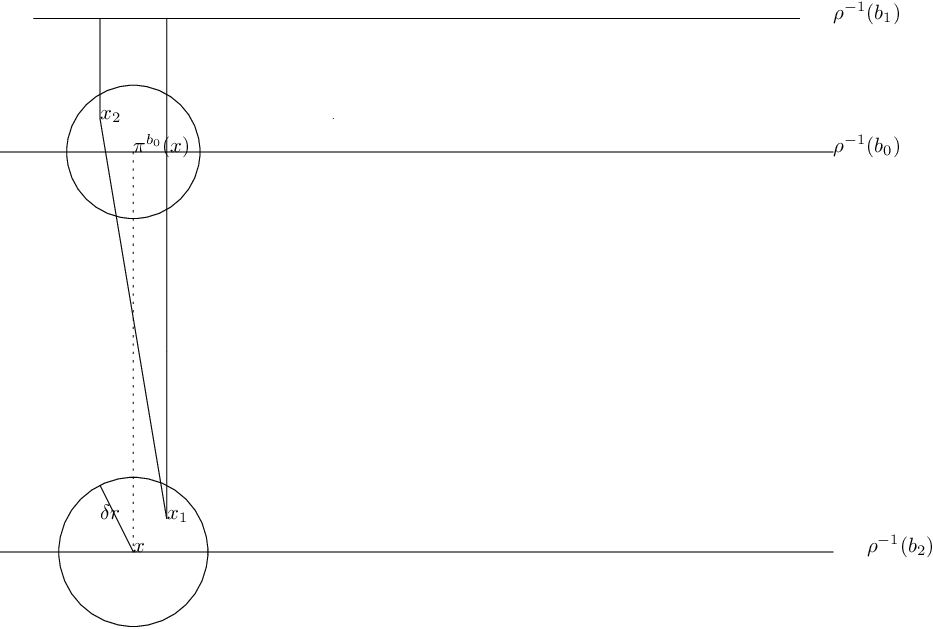}
\caption{Distance estimate of $|x_1x_2|$}
\label{figure: lem Gou-Gu two balls}
\end{center}
\end{figure}

\pf
{Similar to the argument of Lemma \ref{lem property 2 of G-H appr}. 
}
\qed

\subsection{The existence of `good' pair points}

Note in the statement of almost Gou-Gu Theorem (see Theorem \ref{thm almost Gou-Gu}), but we do not know whether
\begin{align}
|\pi^{b_0}(\hat{x})\pi^{b_0}(x)|+ |\pi^{b_0}(\hat{y})\pi^{b_0}(y)|\leq C(n)r\cdot \zeta(\sup_{B_x( 4r)}\sqrt{\frac{|h- \rho|}{r}},  \fint_A H^2,  \fint_A G^2).  \label{proj is really close}
\end{align}
That is why $\hat{x}, \hat{y}$ appear in Theorem \ref{thm almost Gou-Gu}, although we know that the statement is unnatural. This issue can be avoided in local estimate; because we can reduce all discussion to $\hat{x}, \hat{y}$ as in \cite{Xu} and \cite{XZ}, although the presentation is more complicated. 

During the process of gluing the local estimate to obtain the global estimate (we follow the strategy of \cite[Theorem $3.6$]{CC-Ann}, see Section \ref{sec width of almost cylinder}), the projection points problem above seems unavoidable; because we need to do projection repeatedly to get the global distance estimate. 

Based the above reason, we call Theorem \ref{thm almost Gou-Gu} as ``Almost Gou-Gu Theorem".  In the rest of this section,  we show that (\ref{proj is really close}) in fact holds in our context.  Hence we can remove $\hat{x}, \hat{y}$ in the statement of Theorem \ref{thm almost Gou-Gu},  which is Theorem \ref{thm Gou-Gu for cylinder type}.  And we call Theorem \ref{thm Gou-Gu for cylinder type} as ``Stable Gou-Gu Theorem",  which is the main result of this section.

\begin{figure}[H]
\begin{center}
\includegraphics{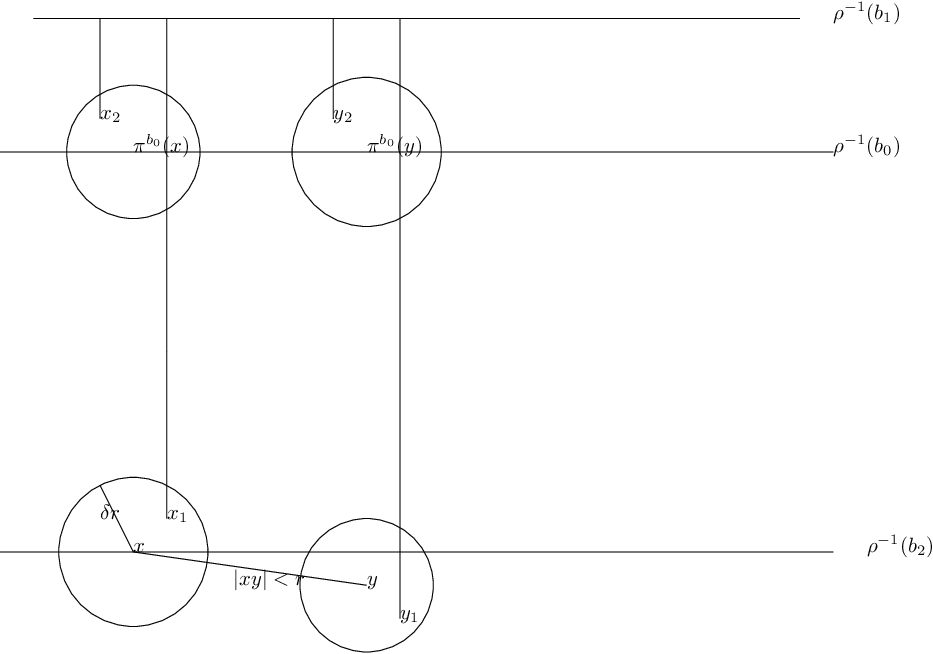}
\caption{Choice of four points $(x_1, x_2, y_1, y_2)$}
\label{figure: prop existence four pts}
\end{center}
\end{figure}

To show (\ref{proj is really close}),  which is the content of Lemma \ref{lem proj dist is controled}; inspired by definition of  the set of `good' points satisfying integral estimate in \cite{CN} (also see \cite{Xu} and \cite{XZ}); we introduce the set of `good' pair points (see Definition \ref{def sets for fours balls around four points} and Figure \ref{figure: prop existence four pts}), whose non-emptiness is proved in Proposition \ref{prop existence of pair of good points}.  

For $\delta> 0, x, y\in B$, let $\Omega_\delta(x, y)= B_x(\delta r)\times B_{\pi^{b_0}(x)}(\delta r)\times B_y(\delta r)\times B_{\pi^{b_0}(y)}(\delta r)$ in the rest of this section.

\begin{definition}\label{def sets for fours balls around four points}
{For $x_1, x_2, y_1, y_2\in B$, we define
\begin{align}
\mathcal{G}(x_1, x_2, y_1, y_2)&\vcentcolon = \sum_{i= 1}^2\int_0^{\rho(x_i)- b_1} G(\xi_{x_i}(t))dt+ \int_0^{\rho(y_i)- b_1} G(\xi_{y_i}(t))dt; \nonumber \\
\mathcal{H}(x_1, x_2, y_1, y_2)&\vcentcolon = \int_0^{\rho(x_1)- b_1}ds(\int_{\overline{\xi_{x_1}(s), \tilde{\xi}_{x_2}(s)}} H)+ \int_0^{\rho(y_1)- b_1}ds(\int_{\overline{\xi_{y_1}(s), \tilde{\xi}_{y_2}(s)}} H) \nonumber \\
&+ \int_0^{\rho(x_1)- b_0}ds(\int_{\overline{\xi_{x_1}(s), \tilde{\xi}_{y_1}(s)}} H). \nonumber
\end{align}
And for $x, y\in B$ we define
\begin{align}
I_{\delta}(x, y)\vcentcolon&= \{(x_1, x_2, y_1, y_2)\in \Omega_\delta(x, y): (\mathcal{G}+ \mathcal{H})(x_1, x_2, y_1, y_2) \leq \delta r\}. \nonumber 
\end{align}
}
\end{definition}

\begin{remark}\label{rem four pts family}
{The philosophy behind this `four points family' is: we not only need to find `good' points near $x, y$ and their projection $\pi(x),  \pi(y)$ with controlled integral,  but also put additional restriction on the integral involving those points respectively.
}
\end{remark}

We recall the following segment inequality due to Cheeger and Colding \cite{CC-Ann}.
\begin{lemma}[Segment Inequality]\label{lem segment ineq}
{Assume $(M^n, g)$ is a complete Riemannian manifold with $Rc\geq 0$, then for any nonnegative function $\Psi$ defined on $B_z( 2r)$, we have
\begin{align}
\int_{B_z( r)\times B_z( r)} \Big(\int_0^{d(y_1, y_2)} \Psi\big(\overline{y_1, y_2}(s)\big) ds\Big) dy_1 dy_2\leq 2^{n+ 2}r\cdot V\big(B_z( r)\big)\cdot \int_{B_z( 2r)} \Psi .\nonumber
\end{align}
}
\end{lemma}\qed

\begin{prop}\label{prop existence of pair of good points}
{There is $C(n)> 0$; such that if $\delta= C(n)\Big\{\big(r^2\fint_A H^2\big)^{\frac{1}{2(2n+ 1)}}+ (\frac{R}{r})(\fint_A G^2)^{\frac{1}{2n}} \Big\}$ where $A= \rho^{-1}(r, 2R- r)$, then for any $|xy|\leq r$, there exists $(\hat{x}, \check{x}, \hat{y}, \check{y})\in I_{\delta}(x, y)$.
}
\end{prop}

\pf
{\textbf{Step (1)}.  Let $\check{I}= B_x(\delta r)\times B_{\pi^{b_0}(x)}(\delta r)\times B_y(\delta r)\times B_{\pi^{b_0}(y)}(\delta r)- I_{\delta}(x, y)$, we only need to show that 
\begin{align}
V(\check{I})< V(B_x(\delta r)\times B_{\pi^{b_0}(x)}(\delta r)\times B_y(\delta r)\times B_{\pi^{b_0}(y)}(\delta r)). \nonumber 
\end{align}

We know
\begin{align}
\delta r\cdot V(\check{I})\leq \int_{\Omega}\mathcal{G}+ \mathcal{H}. \label{G and H integral bound volue of bad pts}
\end{align}
And there is 
\begin{align}
&\int_\Omega dx_1dx_2dy_1dy_2 \int_0^{\rho(x_1)- b_1}G(\xi_{x_1}(t))dt \nonumber \\
&\leq V(B_{\pi^{b_0}(x)}(\delta r)\times B_y(\delta r)\times B_{\pi^{b_0}(y)}(\delta r))\cdot \int_{B_x(\delta r)}dx_1 \int_0^{\rho(x_1)- b_1}G(\xi_{x_1}(t))dt .  \nonumber 
\end{align}

By the mean convex property, we know that the volume element is increasing along $\xi_x(t)$, which is the gradient curve of $-\rho$. Now we get 
\begin{align}
&\int_{B_x(\delta r)}dx_1 \int_0^{\rho(x_1)- b_1}G(\xi_{x_1}(t))dt\leq \int_{b_2- \delta r}^{b_2+ \delta r}ds\int_{\rho^{-1}(s)}dx_1\int_0^{s- b_1}G(\xi_{x_1}(t))dt \nonumber \\
&= \int_{b_2- \delta r}^{b_2+ \delta r}ds \int_0^{s- b_1} dt\int_{\rho^{-1}(s)}G(\xi_{x_1}(t))dx_1 \nonumber \\
&\leq \int_{b_2- \delta r}^{b_2+ \delta r}ds \int_{b_1}^{s} d\tau\int_{\rho^{-1}(\tau)}G(y)dy \leq 2\delta r\cdot \int_{b_1}^{b_2+ 2r} d\tau\int_{\rho^{-1}(\tau)}G(y)dy \nonumber \\
&= C(n)\delta r\cdot \int_{\rho^{-1}(b_1,  b_2+ \delta r)}G(y)dy .\nonumber
\end{align}

Combining the above, we get
\begin{align}
&\frac{\int_{\Omega} dx_1dx_2dy_1dy_2 \int_0^{\rho(x_1)- b_1}G(\xi_{x_1}(t))dt }{V(\Omega)}\leq \frac{C(n)\delta r\cdot \int_{\rho^{-1}(b_1,  b_2+ \delta r)}G(y)dy}{V(B_x(\delta r))}  \nonumber \\
&\leq \frac{C(n)\delta r\cdot \int_{\rho^{-1}(r, 2R- r)}G(y)dy}{V(A)} \cdot \frac{V(A)}{V(B_x(\delta r))} \nonumber \\
&\leq  \frac{C(n)\delta r\cdot \int_{\rho^{-1}(r, 2R- r)}G(y)dy}{V(A)} \cdot \frac{V(B_p(R))}{V(B_x(\delta r))}\nonumber \\
&\leq C(n)(\delta r)\cdot (\frac{R}{\delta r})^n\cdot (\fint_A G^2)^{\frac{1}{2}}.  \label{1st G est}
\end{align}

Similarly we can obtain 
\begin{align}
&\int_\Omega dx_1dx_2dy_1dy_2 \int_0^{\rho(x_2)- b_1}G(\xi_{x_2}(t))dt \nonumber \\
&\leq V(B_{x}(\delta r)\times B_y(\delta r)\times B_{\pi^{b_0}(y)}(\delta r))\cdot C(n)\delta r\int_{\rho^{-1}(b_1, b_0+ \delta r)} G;  \nonumber \\
&\int_\Omega dx_1dx_2dy_1dy_2 \int_0^{\rho(y_1)- b_1}G(\xi_{y_1}(t))dt \nonumber \\
&\leq V(B_{x}(\delta r)\times B_{\pi^{b_0}(x)}(\delta r)\times B_{\pi^{b_0}(y)}(\delta r))\cdot C(n)\delta r\int_{\rho^{-1}(b_1, b_2+ r+ \delta r)} G;  \nonumber \\
&\int_\Omega dx_1dx_2dy_1dy_2 \int_0^{\rho(y_2)- b_1}G(\xi_{y_2}(t))dt \nonumber \\
&\leq V(B_{x}(\delta r)\times B_{\pi^{b_0}(x)}(\delta r)\times B_{y}(\delta r))\cdot C(n)\delta r\int_{\rho^{-1}(b_1, b_0+ \delta r)} G.  \nonumber 
\end{align}

Then we get 
\begin{align}
\frac{\int_{\Omega} \mathcal{G} dx_1dx_2dy_1dy_2 }{V(\Omega)}\leq C(n)(\delta r)\cdot (\frac{R}{\delta r})^n\cdot(\fint_A G^2)^{\frac{1}{2}}.  \label{est of G integral}
\end{align}

\textbf{Step (2)}. Note we have
\begin{align}
&\int_\Omega \int_0^{\rho(x_1)- b_0}dt \int_{\overline{\xi_{x_1}(t), \tilde{\xi}_{y_1}(t)}} H \nonumber \\
&\leq V(B_{\pi(x)}(\delta r)\times B_{\pi(y)}(\delta r))\int_{B_{x}(\delta r)\times B_{y}(\delta r)}dx_1dy_1 \int_0^{\rho(x_1)- b_0}dt \int_{\overline{\xi_{x_1}(t), \tilde{\xi}_{y_1}(t)}} H.\nonumber
\end{align}

From the Co-area formula and volume comparison between level sets of $\rho$,  we have
\begin{align}
&\int_{B_{x}(\delta r)\times B_{y}(\delta r)}dx_1dy_1 \int_0^{\rho(x_1)- b_0}dt \int_{\overline{\xi_{x_1}(t), \tilde{\xi}_{y_1}(t)}} H \nonumber\\
&\leq \int_{b_2- \delta r}^{b_2+ \delta r} dt_1\int_{\rho(y)- \delta r}^{\rho(y)+ \delta r}  dt_2 \int_{\rho^{-1}(t_1)\times \rho^{-1}(t_2)\cap (B_x(\delta r)\times B_y(\delta r))} dx_1 dy_2 \int_0^{t_1- b_0}dt\Big(\int_{\overline{\xi_{x_1}(t), \tilde{\xi}_{y_1}(t)}} H\Big) ds \nonumber \\
&\leq \int_{0}^{b_2- b_0+ \delta r} dt  \int_{B_x( 4r)\times B_x( 4r)} d\tilde{x}d\tilde{y} \Big(\int_{\overline{\tilde{x}, \tilde{y}}} H\Big) \nonumber \\
&\leq  C(n)r^2V(B_x(4r))\int_{B_x(8r)} H .\nonumber 
\end{align}

Now we get
\begin{align}
&\frac{\int_\Omega \int_0^{\rho(x_1)- b_0}dt \int_{\overline{\xi_{x_1}(t), \tilde{\xi}_{y_1}(t)}} H}{V(\Omega)}\leq \frac{C(n)r^2V(B_x(4r))\int_{B_x(8r)} H}{V(B_x(\delta r))\cdot V(B_y(\delta r))} \leq C(n)\delta^{-2n}r^{2} (\fint_A H^2)^{\frac{1}{2}}.  \nonumber 
\end{align}

Similarly we can obtain
\begin{align}
&\frac{\int_\Omega \Big( \int_0^{\rho(x_1)- b_1}dt \int_{\overline{\xi_{x_1}(t), \tilde{\xi}_{x_2}(t)}} H+ \int_0^{\rho(y_1)- b_1}dt \int_{\overline{\xi_{y_1}(t), \tilde{\xi}_{y_2}(t)}} H\Big) }{V(\Omega)} \leq C(n)\delta^{-2n}r^{2} (\fint_A H^2)^{\frac{1}{2}}.  \nonumber 
\end{align}

Therefore we have
\begin{align}
\frac{\int_{\Omega} \mathcal{H} dx_1dx_2dy_1dy_2 }{V(\Omega)}\leq C(n)\delta^{-2n}r^{2} (\fint_A H^2)^{\frac{1}{2}}. \label{est of H integral}
\end{align}

Combining (\ref{G and H integral bound volue of bad pts}),  (\ref{est of G integral}) with (\ref{est of H integral}),  we get
\begin{align}
\frac{V(\check{I})}{V(\Omega)}&\leq \frac{\int_\Omega \mathcal{G}+ \mathcal{H}}{\delta r V(\Omega)} \leq C(n)\Big\{\delta^{-2n- 1}r(\fint_A H^2)^{\frac{1}{2}}+ \big(\frac{R}{\delta r}\big)^n(\fint_A G^2)^{\frac{1}{2}}\Big\}.  \nonumber 
\end{align}

Then there is $C(n)> 0$ such that if 
\begin{align}
\delta= C(n)\Big\{\big(r^2\fint_A H^2\big)^{\frac{1}{2(2n+ 1)}}+ (\frac{R}{r})(\fint_A G^2)^{\frac{1}{2n}} \Big\}, \nonumber 
\end{align}
then $\frac{V(\check{I})}{V(\Omega)}< 1$ and the conclusion follows.
}
\qed

\subsection{The local stable Gou-Gu theorem}

Therefore Theorem \ref{thm Gou-Gu for cylinder type} can be proved by using Lemma \ref{lem proj dist is controled} and Almost Gou-Gu Theorem, and we provide a unified argument to obtain Theorem \ref{thm Gou-Gu for cylinder type} for self-completeness. 

In this subsection we assume $\delta= C(n)\Big\{\big(r^2\fint_A H^2\big)^{\frac{1}{2(2n+ 1)}}+ (\frac{R}{r})(\fint_A G^2)^{\frac{1}{2n}} \Big\}$,  where $C(n)$ is as in Proposition \ref{prop existence of pair of good points}.

\begin{lemma}\label{lem dist for good pair of x and y}
{For $(\hat{x}, \check{x}, \hat{y}, \check{y})\in I_{\delta}(x, y)$, we have 
\begin{align}
\big||\hat{y}\hat{x}|- |\mathscr{G}^{b_0}(\hat{y}) \mathscr{G}^{b_0}(\hat{x})|\big|\leq C(n)r\cdot (\delta^{\frac{1}{2}}+ \sup_{B_x(8r)}\frac{|h- \rho|}{r}). \nonumber 
\end{align}
}
\end{lemma}

\begin{figure}[H]
\begin{center}
\includegraphics{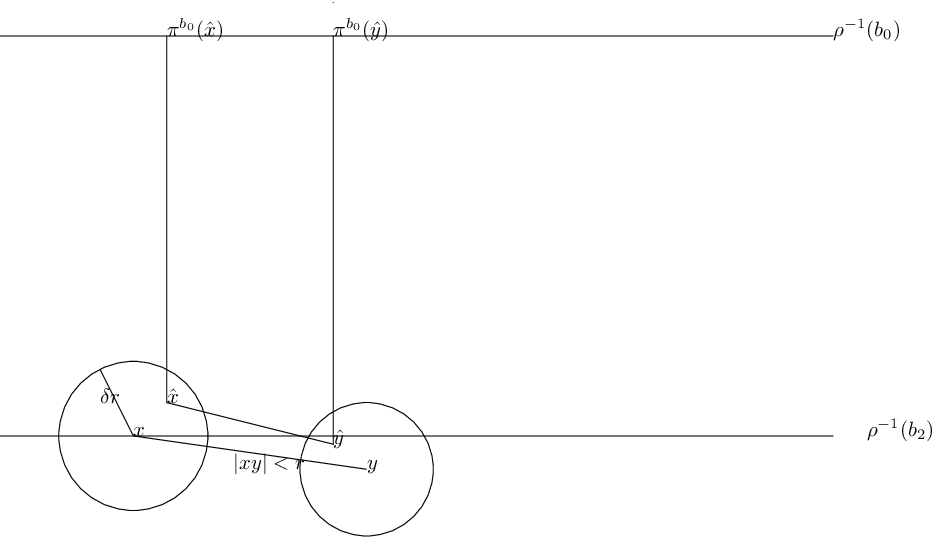}
\caption{Distance estimate of good pair points $|\hat{y}\hat{x}|$}
\label{figure: lem dist for good x y}
\end{center}
\end{figure}

\pf
{From Lemma \ref{lem property 2 of G-H appr}, we know that
\begin{align}
&||\hat{y}\hat{x}|^2- |\mathscr{G}^{b_0}(\hat{y}) \mathscr{G}^{b_0}(\hat{x})|^2| \nonumber \\
&\leq 16r\cdot(\int_0^{\rho(\hat{y})- b_0} \big|\nabla (\rho- h)\big|(\xi_{\hat{y}}(t))dt+ \int_0^{\rho(\hat{x})- b_0}\big|\nabla (\rho- h)\big|(\xi_{\hat{x}}(t))dt) \nonumber \\
&\quad  +4|\rho(\hat{x})- \rho(\hat{y})|\cdot \sup_{B_x( 4r)}|h- \rho|+ 64r\cdot \int_0^{\rho(\hat{x})- b_0} dt\int_{0}^{l_t} |\nabla^2 h|(\tau_t(s)) ds\nonumber \\
&\leq C(n)\cdot (\delta r^2+ |\rho(\hat{x})- \rho(\hat{y})|\cdot \sup_{B_x( 4r)}|h- \rho|). \label{the result of Gou-Gu}
\end{align}

Now we note $|\nabla \rho|= 1$, then from the above we get 
\begin{align}
(I)\vcentcolon&= ||\hat{y}\hat{x}|- |\mathscr{G}^{b_0}(\hat{y}) \mathscr{G}^{b_0}(\hat{x})||= \frac{||\hat{y}\hat{x}|^2- |\mathscr{G}^{b_0}(\hat{y}) \mathscr{G}^{b_0}(\hat{x})|^2|}{|\hat{y}\hat{x}|+ |\mathscr{G}^{b_0}(\hat{y}) \mathscr{G}^{b_0}(\hat{x})|} \nonumber \\
&\leq C(n)\cdot \Big\{\frac{\delta r^2}{|\hat{y}\hat{x}|+ |\mathscr{G}^{b_0}(\hat{y}) \mathscr{G}^{b_0}(\hat{x})|}+ \sup_{B_x( 4r)}|h- \rho|)\Big\}. \label{before seperate cases}
\end{align}

There are two cases to be dealt with separately as follows:
\begin{enumerate}
\item[Case (1)]: If $|\hat{x}\hat{y}|\geq \tau r$, where $\tau> 0$ is to be determined later; from (\ref{before seperate cases}) we have 
\begin{align}
(I)\leq C(n)\cdot \Big\{\frac{\delta r}{\tau}+ \sup_{B_x( 4r)}|h- \rho|)\Big\}. \label{case 1 est}
\end{align}
\item[Case (2)]: If $|\hat{x}\hat{y}|\leq \tau r$, from (\ref{the result of Gou-Gu}) we have 
\begin{align}
|\mathscr{G}^{b_0}(\hat{y})\mathscr{G}^{b_0}(\hat{x})|^2&\leq |\hat{x}\hat{y}|^2+ C(n)\cdot (\delta r^2+ |\rho(\hat{x})- \rho(\hat{y})|\cdot \sup_{B_x( 4r)}|h- \rho|) \nonumber \\
&\leq (\tau r)^2+ C(n)\cdot (\delta r^2+ \tau r\cdot \sup_{B_x( 4r)}|h- \rho|), \nonumber 
\end{align}
which implies 
\begin{align}
|\mathscr{G}^{b_0}(\hat{y})\mathscr{G}^{b_0}(\hat{x})|&\leq C(n)r\cdot \Big\{\tau + \delta^{\frac{1}{2}}+ \sqrt{\tau\cdot \sup_{B_x( 4r)}\frac{|h- \rho|}{r}}\Big\} \nonumber \\
&\leq C(n)r\cdot \Big\{\tau + \delta^{\frac{1}{2}}+  \sup_{B_x( 4r)}\frac{|h- \rho|}{r}\Big\} . \nonumber 
\end{align}

Finally we get
\begin{align}
(I)\leq |\hat{y}\hat{x}|+ |\mathscr{G}^{b_0}(\hat{y})\mathscr{G}^{b_0}(\hat{x})|\leq C(n)r\cdot \Big\{\tau + \delta^{\frac{1}{2}}+  \sup_{B_x( 4r)}\frac{|h- \rho|}{r}\Big\}. \nonumber 
\end{align}
\end{enumerate}

By Case (1) and (2) above, we always have
\begin{align}
(I)\leq C(n)r\cdot \Big\{\tau + \frac{\delta}{\tau}+ \delta^{\frac{1}{2}}+  \sup_{B_x( 4r)}\frac{|h- \rho|}{r}\Big\}. \nonumber 
\end{align}

Let $\tau= \frac{\delta}{\tau}$, then $\tau= \delta^{\frac{1}{2}}$. Plug the value of $\tau$ into the above inequality, we obtain
\begin{align}
(I)\leq C(n)r\cdot \Big\{\delta^{\frac{1}{2}}+  \sup_{B_x( 4r)}\frac{|h- \rho|}{r}\Big\}. \nonumber 
\end{align}
}
\qed

\begin{lemma}\label{lem proj dist is controled}
{For $(\hat{x}, \check{x}, \hat{y}, \check{y})\in I_{\delta}(x, y)$, we have
\begin{align}
|\mathscr{G}^{b_0}(x)\mathscr{G}^{b_0}(\hat{x})|+ |\mathscr{G}^{b_0}(y)\mathscr{G}^{b_0}(\hat{y})|\leq C(n)r\cdot \Big\{\sup_{B_x( 4r)}\sqrt{\frac{|h- \rho|}{r}}+ \sqrt{\delta}\Big\} .  \nonumber 
\end{align}
}
\end{lemma}

\begin{figure}[H]
\begin{center}
\includegraphics{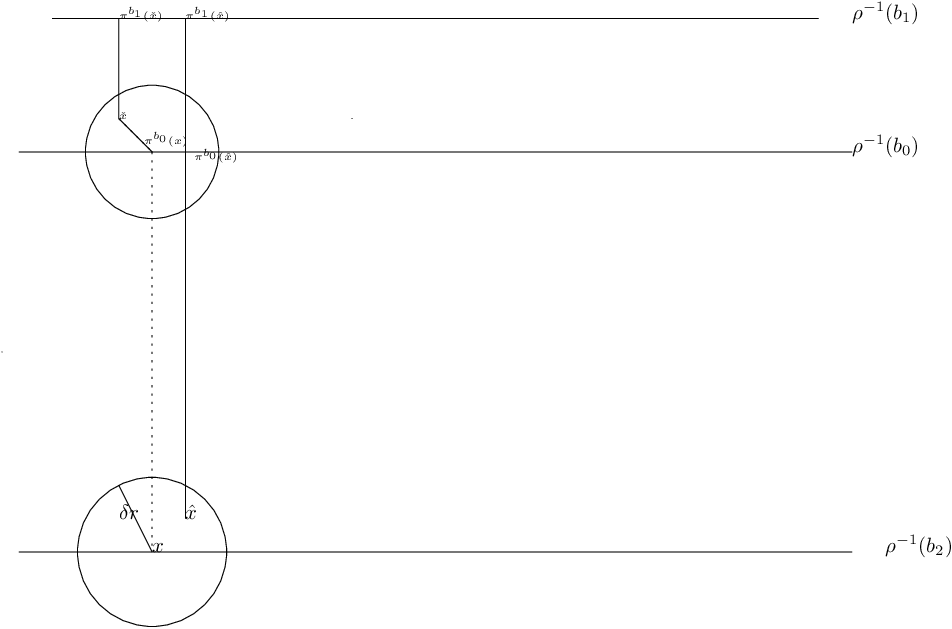}
\caption{Distance estimate of $|\pi^{b_0}(x)\pi^{b_0}(\hat{x})|$}
\label{figure: lem proj est}
\end{center}
\end{figure}

\pf
{\textbf{Step (1)}. From Lemma \ref{lem Gou-Gu for two balls} and the definition \ref{def sets for fours balls around four points}, we have
\begin{align}
&|\pi^{b_1}(\hat{x})\pi^{b_1}(\check{x})|^2\leq |\hat{x}\check{x}|^2- |\rho(\hat{x})- \rho(\check{x})|^2+ 16r\cdot(\int_0^{\rho(\check{x})- b_1} G(\xi_{\check{x}}(t))dt+ \int_0^{\rho(\hat{x})- b_1}G(\xi_{\hat{x}}(t))dt) \nonumber \\
&\quad  +4|\rho(\hat{x})- \rho(\check{x})|\cdot \sup_{B_x( 4r)}|h- \rho|+ 64r\cdot \int_0^{\rho(\hat{x})- b_1} dt\int_{0}^{l_t} H(\tau_t(s)) ds \nonumber \\
&\leq (|x\pi^{b_0}(x)|+ 2\delta r)^2- (\rho(x)- 2\delta r- b_0)^2+ C(n)r^2\cdot \Big\{\sup_{B_x( 4r)}\frac{|h- \rho|}{r}+ \delta \Big\} \nonumber \\
&\leq C(n)\delta r|\rho(x)- b_0|+ C(n)r^2\cdot \Big\{\sup_{B_x( 4r)}\frac{|h- \rho|}{r}+ \delta \Big\}  \nonumber \\
&\leq  C(n)r^2\cdot \Big\{\sup_{B_x( 4r)}\frac{|h- \rho|}{r}+ \delta \Big\} . \nonumber 
\end{align}

Therefore we have
\begin{align}
|\pi^{b_1}(\hat{x})\pi^{b_1}(\check{x})|\leq C(n)r\cdot \Big\{\sup_{B_x( 4r)}\sqrt{\frac{|h- \rho|}{r}}+  \sqrt{\delta}\Big\}. \label{the proj of hat x and check x are close}
\end{align}

\textbf{Step (2)}. By metric triangle inequality and (\ref{the proj of hat x and check x are close}), we get 
\begin{align}
|\pi^{b_0}(x)\pi^{b_0}(\hat{x})|&\leq |\pi^{b_0}(x)\check{x}|+ |\check{x}\pi^{b_1}(\check{x})|+ |\pi^{b_1}(\hat{x})\pi^{b_1}(\check{x})|+ |\pi^{b_1}(\hat{x})\pi^{b_0}(\hat{x})| \nonumber \\
&\leq 2\delta r+ 2|b_1- b_0|+ C(n)r\cdot \Big\{\sup_{B_x( 4r)}\sqrt{\frac{|h- \rho|}{r}}+ \sqrt{\delta}+ \sqrt{\delta}\Big\} \nonumber \\
&\leq C(n)r\cdot \Big\{\sup_{B_x( 4r)}\sqrt{\frac{|h- \rho|}{r}}+ \sqrt{\delta} \Big\}. \nonumber 
\end{align}

Then note $|\rho(x)- \rho(\hat{x})| \leq |x\hat{x}|\leq \delta r$, we get 
\begin{align}
|\mathscr{G}^{b_0}(x)\mathscr{G}^{b_0}(\hat{x})|&\leq |\pi^{b_0}(x)\pi^{b_0}(\hat{x})|+ |\rho(x)- \rho(\hat{x})| \leq C(n)r\cdot \Big\{\sup_{B_x( 4r)}\sqrt{\frac{|h- \rho|}{r}}+ \sqrt{\delta}\Big\}. \nonumber 
\end{align}
}
\qed

\begin{theorem}\label{thm Gou-Gu for cylinder type}
{There is $C(n)> 0$, such that for any $|xy|\leq r$ where $r<< 1$, we have
\begin{align}
&\Big|d(x, y)- |\mathscr{G}^{b_0}(x)\mathscr{G}^{b_0}(y)|\Big| \nonumber \\
&\leq C(n)r\cdot \Big\{\sup_{B_x( 4r)}\sqrt{\frac{|h- \rho|}{r}}+  \big(r^2\fint_A H^2\big)^{\frac{1}{4(2n+ 1)}}+ (\frac{R}{r})^{\frac{1}{2}}(\fint_A G^2)^{\frac{1}{4n}}\Big\}. \nonumber 
\end{align}
}
\end{theorem}

\pf
{From Proposition \ref{prop existence of pair of good points} we find $(\hat{x}, \check{x}, \hat{y}, \check{y})\in I_{\delta}(x, y)$.

Now from Lemma \ref{lem dist for good pair of x and y} (distance estimate between good pair of points) and Lemma \ref{lem proj dist is controled} (`model' distance estimate between original points and good points), we have 
\begin{align}
&\quad \Big||xy|- |\mathscr{G}^{b_0}(x)\mathscr{G}^{b_0}(y)|\Big| \nonumber \\
&\leq \Big||xy|- |\hat{x}\hat{y}|\Big|+ \Big||\hat{x}\hat{y}|- |\mathscr{G}^{b_0}(\hat{x})\mathscr{G}^{b_0}(\hat{y})|\Big|+ \Big||\mathscr{G}^{b_0}(\hat{x})\mathscr{G}^{b_0}(\hat{y})|- |\mathscr{G}^{b_0}(x)\mathscr{G}^{b_0}(y)|\Big|  \nonumber \\
&\leq |x\hat{x}|+ |y\hat{y}|+ C(n)r\cdot (\delta^{\frac{1}{2}}+ \sup_{B_x(8r)}\frac{|h- \rho|}{r})+ |\mathscr{G}^{b_0}(x)\mathscr{G}^{b_0}(\hat{x})|+ |\mathscr{G}^{b_0}(y)\mathscr{G}^{b_0}(\hat{y})|\nonumber \\
&\leq 2\delta r+ C(n)r\cdot (\delta^{\frac{1}{2}}+ \sup_{B_x(8r)}\frac{|h- \rho|}{r})+ C(n)r\cdot \Big\{\sup_{B_x( 4r)}\sqrt{\frac{|h- \rho|}{r}}+ \sqrt{\delta}\Big\} \nonumber \\
&\leq C(n)r\cdot \Big\{\sup_{B_x( 4r)}\sqrt{\frac{|h- \rho|}{r}}+  \sqrt{\delta}\Big\}\nonumber \\
&\leq C(n)r\cdot \Big\{\sup_{B_x( 4r)}\sqrt{\frac{|h- \rho|}{r}}+  \big(r^2\fint_A H^2\big)^{\frac{1}{4(2n+ 1)}}+ (\frac{R}{r})^{\frac{1}{2}}(\fint_A G^2)^{\frac{1}{4n}} \Big\}. \nonumber 
\end{align}
}
\qed

\begin{prop}\label{prop Gou-Gu with spectral gap}
{There is $C(n)> 0$, such that for any $x, y$ with $|xy|< bR$ where $b<< 1$, we have 
\begin{align}
\Big|d(x, y)- |\mathscr{G}^{b_0}(x)\mathscr{G}^{b_0}(y)|\Big|\leq C(n)(bR)\Big\{b^{-1+ \frac{1}{2n+ 4}}(\epsilon R^2)^{\frac{1}{8(2n+ 1)}}\Big\} .\nonumber 
\end{align}
}
\end{prop}

\pf
{It directly follows from Theorem \ref{thm Gou-Gu for cylinder type} and Lemma \ref{lem grad and Hessian est of h}.  
}
\qed

\section{The width estimate of almost cylinder}\label{sec width of almost cylinder}

In this paper, to get an explicit width estimate of geodesic ball, we need to obtain global distance estimate. Here the global distance estimate means that the segment between $x, y$ possibly touch the boundary of geodesic ball (if the geodesic ball belongs to some complete Riemannian manifold, the segment between $x, y$ possibly run out of the geodesic ball). Anyway, we need to drop the closeness assumption of $x, y$ in the local distance estimate above. The method of getting the global distance estimate from the local distance estimate, is discussed in \cite{CC-Ann} firstly. 

In this section, we obtain the global distance estimate, using local stable Gou-Gu theorem and the property that the gradient curves of distance functions are almost extendable. The philosophy of our argument is close to the strategy of \cite[Theorem $3.6$]{CC-Ann}. Because we focus on the concrete case and are interested in obtaining the quantitative estimate of diameter of the level sets, we give the detailed discussion in our context.

Let $I_k= [kbR, (k+ 1)bR]$, where $k\in \mathbb{Z}^+$. We define $\pi^{(k)}(x)= \xi_x\cap \rho^{-1}(kbR)$ and 
\begin{align}
\mathscr{G}^{(k)}(x)= (\pi^{(k)}(x), \rho(x)). \nonumber 
\end{align}

The main result of this section is the following Theorem. 
\begin{theorem} 
{Assume there is $\delta_1, \delta_2> 0$,  such that, for any $k\in \mathbb{Z}^+$ and two points $x, y$ satisfying
\begin{align}
y\in B(x, bR), \quad \quad \rho(x)- kbR\leq \frac{3}{2}bR, \quad \quad \min\{\rho(x)- kbR, \rho(y)- kbR\}\geq \frac{1}{2}bR, \nonumber 
\end{align}
we have 
\begin{align}
\Big||xy|- |\mathscr{G}^{(k)}(x)\mathscr{G}^{(k)}(y)|\Big|\leq \delta_1 bR. \label{assumed Gou-Gu}
\end{align}
Further assume that for any $x\in \rho^{-1}(bR)$ and $t\in [b, 2-b]$, there is $z\in \rho^{-1}(tR)$ such that 
\begin{align}
d(x, z)- (t- b)R\leq \delta_2 R. \label{geodesic line extendable}
\end{align}
Then $\displaystyle \mathcal{W}(B_p(R))\leq 10^7 R\Big\{b+ \sqrt{\frac{\delta_1}{b}+ \delta_2}\Big\}$.
}
\end{theorem}

\subsection{The estimate of the outer level sets' diameter}

\begin{lemma}\label{lem almost proj is close to true proj}
{If $B_p(R)$ is a metric ball satisfying (\ref{assumed Gou-Gu}), then for any $y\in \rho^{-1}(8bR), t\in [8b, 2-b], w\in \rho^{-1}(tR)$ satisfying
\begin{align}
d(y, w)- (t- 8b)R\leq \delta_2 R; \nonumber 
\end{align}
we have
\begin{align}
&\{z_i\}_{i= 0}^m \subseteq \rho^{-1}(bR), \quad  z_m= \pi^{bR}(w),  \quad  z_0= \pi^{bR}(y),  \nonumber\\
&\sum_{i= 1}^m |z_iz_{i- 1}|\leq 20R(b+ \sqrt{\frac{\delta_1}{b}+ \delta_2})\sqrt{t- 7b}; \nonumber \\
&|\pi^{8bR}(w)y|\leq 20R\Big\{b+ \sqrt{\frac{\delta_1}{b}+ \delta_2}\sqrt{t- 7b}\Big\}; \nonumber 
\end{align}
where $m= \Big[\frac{t- 8b}{b}\Big]+ 1$ and $\Big[\frac{t- 8b}{b}\Big]$ is the biggest integer no more than $\frac{t- 8b}{b}$.
}
\end{lemma}

\pf 
{\textbf{Step (1)}. Let $y_0= y$, choose $\{y_i\}_{i= 1}^{m}$ from geodesic segment $\overline{w, y_0}$, where 
\begin{align}
y_m= w, \quad \quad |y_iy_0|= \frac{i}{m}|wy_0|.\nonumber 
\end{align}

From the assumption, we know that
\begin{align}
(t- 8b)R\leq |wy|\leq (t- 8b+ \delta_2)R. \label{wy two bounds}
\end{align}

Using (\ref{wy two bounds}) we have
\begin{align}
\rho(y_i)&\geq \rho(w)- |wy_i|\geq tR- (1- \frac{i}{m})|wy|\geq 8bR+ \frac{i}{m}(t- 8b)R- b(1- \frac{i}{m})R \nonumber \\
&\geq 8bR+ \frac{i}{m}|wy|- b R . \label{lower bound 2}
\end{align}

Also we have the upper bound
\begin{align}
\rho(y_i)\leq \rho(y_0)+ |y_iy_0|\leq 8bR+ \frac{i}{m}|wy|. \label{upper bound 2}
\end{align}

\textbf{Step (2)}. For each $1\leq i\leq m$, we define
\begin{align}
\tau(i)= \max\{j\in \overline{\mathbb{Z}^+}: \min\{\rho(y_i)- jbR, \rho(y_{i- 1})- jbR\}\geq \frac{1}{2}bR\}. \label{def of tau i}
\end{align}

Note the above choice of $\tau(i)$ also implies
\begin{align}
\min\{\rho(y_i)-\tau(i)bR, \rho(y_{i- 1})- \tau(i)bR\}\leq \frac{3}{2}bR. \label{another bound of tau i}
\end{align}

From (\ref{assumed Gou-Gu}) and the choice of $\tau(i)$, we get 
\begin{align}
\Big||y_{i- 1}y_i|- |\mathscr{G}^{(\tau(i))}(y_{i- 1})\mathscr{G}^{(\tau(i))}(y_i)|\Big|\leq \delta_1 bR, \label{Gou-Gu applies on yi}
\end{align}

By (\ref{def of tau i})  and (\ref{upper bound 2}), we obtain
\begin{align}
(\tau(i)+ \frac{1}{2})bR\leq 8bR+ \frac{i}{m}|wy|. \nonumber 
\end{align}
Using the definition of $m$, the above inequality implies
\begin{align}
\tau(i)\leq i+ 9. \label{upper bound of tau i-final}
\end{align}

Similarly, from (\ref{another bound of tau i}) and (\ref{lower bound 2}) we get
\begin{align}
(\tau(i)+ \frac{3}{2})bR\geq 8bR+ \frac{i- 1}{m}|wy|- bR. \nonumber 
\end{align}
Again using the definition of $m$, we further obtain
\begin{align}
\tau_i\geq i+ 3. \label{lower bound of tau i-final}
\end{align}

We define 
\begin{equation}\nonumber 
y_i^{\tau(i), l}\vcentcolon= \left\{
\begin{array}{rl}
&\pi^{(l)}\circ \pi^{(l+ 1)}\cdots \pi^{(\tau(i))}(y_i) , \quad \quad \quad \quad \quad l\leq \tau(i),  \\
&\pi^{(l)}\circ \pi^{(l- 1)}\cdots \pi^{(\tau(i))}(y_i),  \quad \quad \quad \quad \quad l\geq \tau(i).  
\end{array} \right.
\end{equation}
From (\ref{assumed Gou-Gu}) we also have
\begin{align}
\Big||y_{i- 1}^{\tau(i), \tau(i)} y_i^{\tau(i), \tau(i)}|- |y_{i- 1}^{\tau(i), 1} y_i^{\tau(i), 1}|\Big|\leq (\tau(i)- 1)\cdot \delta_1 bR\leq (i+ 8)\delta_1 bR. \label{parallel Gou-Gu application}
\end{align}

\begin{figure}[H]
\begin{center}
\includegraphics{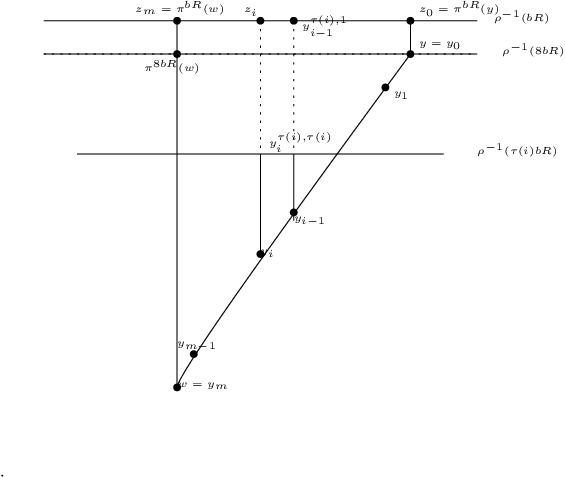}
\caption{Multiple projections of $y_i$}
\label{figure: prop general-1}
\end{center}
\end{figure}

By (\ref{Gou-Gu applies on yi}) and (\ref{parallel Gou-Gu application}), we obtain
\begin{align}
&\quad \Big||y_{i- 1}y_i|^2- |\rho(y_{i- 1})- \rho(y_i)|^2- |y_{i- 1}^{\tau(i), 1} y_i^{\tau(i), 1}|^2\Big| \nonumber \\
&\leq \Big||y_{i- 1}y_i|^2- |\rho(y_{i- 1})- \rho(y_i)|^2- |y_{i- 1}^{\tau(i), \tau(i)} y_i^{\tau(i), \tau(i)}|^2\Big|+ \Big||y_{i- 1}^{\tau(i), \tau(i)} y_i^{\tau(i), \tau(i)}|^2- |y_{i- 1}^{\tau(i), 1} y_i^{\tau(i), 1}|^2\Big| \nonumber \\
&\leq 4(bR)^2\delta_1+ 4(i+ 8)\delta_1 (bR)^2. \nonumber
\end{align}

\textbf{Step (3)}. Take the sum with respect to $i$, note (\ref{upper bound of tau i-final}) and (\ref{lower bound of tau i-final}),  we have
\begin{align}
&\quad \sum_{i= 1}^{m}|y_{i- 1}^{\tau(i), 1} y_i^{\tau(i), 1}|^2\leq \sum_{i= 1}^{m}|y_{i- 1}y_i|^2- \sum_{i= 1}^{m}|\rho(y_{i- 1})- \rho(y_i)|^2 + 4\delta_1(bR)^2 \sum_{i= 1}^m (i+ 9) \nonumber \\
&\leq \frac{|wy_0|^2}{m}- \frac{(\sum\limits_{i= 1}^{m}|\rho(y_{i- 1})- \rho(y_i)|)^2}{m}+ 36(bR)^2\delta_1 m^2 \nonumber \\
&\leq \frac{1}{m}((t- 8b+ \delta_2)R)^2- \frac{1}{m}|\rho(y_0)- \rho(w)|^2+ 200\delta_1R^2\leq 200R^2(\delta_2\cdot b + \delta_1) . \nonumber 
\end{align}

On the other hand, by Cauchy-Schwarz inequality we have 
\begin{align}
\sum_{i= 1}^{m}|y_{i- 1}^{\tau(i), 1} y_i^{\tau(i), 1}|^2\geq \frac{1}{m}\Big(\sum_{i= 1}^{m}|y_{i- 1}^{\tau(i), 1} y_i^{\tau(i), 1}|\Big)^2\geq \frac{b}{t- 7b}|\pi^{bR}(y)w^{m, 1}|^2. \nonumber 
\end{align}

Note $w^{m, 1}= \pi^{bR}(w)$ and $y_0= y$, combining the above yields
\begin{align}
|\pi^{bR}(y)\pi^{bR}(w)|\leq 20R\sqrt{\frac{\delta_1}{b}+ \delta_2}\sqrt{t- 7b}. \nonumber 
\end{align}

Therefore
\begin{align}
|\pi^{8bR}(w)y|& \leq |\pi^{bR}(y)\pi^{bR}(w)|+ |\pi^{bR}(y)y|+ |\pi^{8bR}(w)\pi^{bR}(w)|  \nonumber \\
&\leq 20R\Big\{b+ \sqrt{\frac{\delta_1}{b}+ \delta_2}\sqrt{t- 7b}\Big\}.  \nonumber 
\end{align}

Let $z_i= \pi^{bR}(y_i)$, the conclusion follows from the above.
}
\qed

\begin{prop}\label{prop outer level set diam upper bound}
{If (\ref{assumed Gou-Gu}) holds, then we have
\begin{align}
\sup_{x, y\in \rho^{-1}(8bR)}|xy|\leq 80R\Big\{b+ \sqrt{\frac{\delta_1}{b}+ \delta_2}\Big\}. \nonumber 
\end{align}
}
\end{prop}

\pf
{In Lemma \ref{lem almost proj is close to true proj}, let $w= p$ there, then $t= 1$. For any $x_0, y_0\in \rho^{-1}(8bR)$, we obtain
\begin{align}
|x_0y_0|&\leq |x_0\pi^{8bR}(p)|+ |y_0\pi^{8bR}(p)|\leq 40R\Big\{b+ \sqrt{\frac{\delta_1}{b}+ \delta_2}\Big\} \nonumber \\
&\leq 80R\Big\{b+ \sqrt{\frac{\delta_1}{b}+ \delta_2}\Big\}. \nonumber 
\end{align}
}
\qed

\subsection{The width estimate of metric ball}

Recall the width of a metric ball is defined as follows:
\begin{align}
\mathcal{W}(B_p(R))= \inf_{f\in \mathrm{Lip}(B_p(R))}\max_{t\in \mathbb{R}}\mathrm{diam}(f^{-1}(t)).  \nonumber 
\end{align}

\begin{theorem}\label{thm diam of the level sets}
{If $B_p(R)$ is a metric ball satisfying (\ref{assumed Gou-Gu}) and (\ref{geodesic line extendable}), then 
\begin{align}
\mathcal{W}(B_p(R))\leq 10^7 R\Big\{b+ \sqrt{\frac{\delta_1}{b}+ \delta_2}\Big\}.  \nonumber 
\end{align}
}
\end{theorem}

\pf
{\textbf{Step (1)}. We firstly consider any $x, y\in \rho^{-1}(tR)$ where $t\in [8b, 2- b]$. Applying Lemma \ref{lem almost proj is close to true proj} on $\pi^{8bR}(x), p$ and $\pi^{8bR}(y), p$ respectively, we can find 
\begin{align}
\{z_i\}_{i= 0}^m \subseteq \rho^{-1}(bR), \quad z_m= \pi^{bR}(p),  \quad z_0= \pi^{bR}(x),   \quad \sum_{i= 1}^m |z_iz_{i- 1}|\leq 20R(b+ \sqrt{\frac{\delta_1 }{b}+ \delta_2})\sqrt{1- 7b}; \nonumber \\
\{z_i\}_{i= m}^{2m} \subseteq \rho^{-1}(bR), \quad  z_m= \pi^{bR}(p),  \quad  z_{2m}= \pi^{bR}(y),   \quad \sum_{i= m}^{2m} |z_iz_{i- 1}|\leq 20R(b+ \sqrt{\frac{\delta_1 }{b}+ \delta_2})\sqrt{1- 7b}. \nonumber 
\end{align}

\begin{figure}[H]
\begin{center}
\includegraphics{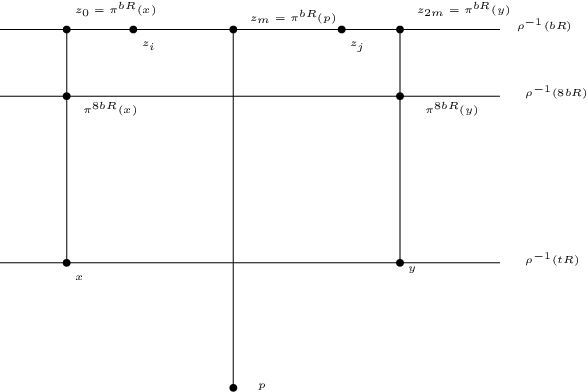}
\caption{The choice of $z_i$}
\label{figure: thm diam of level sets}
\end{center}
\end{figure}

We define $\check{x}= \pi^{bR}(x),  \check{y}= \pi^{bR}(y)$. Therefore we can find $\{v_i\}_{i= 0}^{400}\subseteq \{z_i\}_{i= 0}^{2m}\subseteq \rho^{-1}(bR)$ such that 
\begin{align}
&v_0= \check{x},  \quad \quad v_{400}= \check{y},  \quad \quad \sum_{i= 1}^{400} |v_iv_{i- 1}|\leq 40R(b+ \sqrt{\frac{\delta_1 }{b}+ \delta_2})\sqrt{t- 7b}; \label{sum of zi zi-1} \\
&|v_iv_{i- 1}|\leq 5^{-1}R(b+ \sqrt{\frac{\delta_1 }{b}+ \delta_2})\sqrt{t- 7b} .  \nonumber 
\end{align}

For each $v_i\in \rho^{-1}(bR)$,  by (\ref{geodesic line extendable}) we can find $w_i\in \rho^{-1}(tR)$ such that
\begin{align}
w_0= x,  \quad \quad w_{400}= y,  \quad \quad \quad |w_iv_i|- (t- b)R\leq b R.  \label{almost proj ineq}
\end{align}

Let $\hat{w}_i= \pi^{bR}(w_i)$,  from (\ref{almost proj ineq}), we apply Lemma \ref{lem almost proj is close to true proj} on $w_i, v_i$ (replacing $b$ in Lemma \ref{lem almost proj is close to true proj} by $8^{-1}b$), then
\begin{align}
|\hat{w}_iv_i|\leq 20R\Big\{8^{-1}b+ (b+ \sqrt{\frac{\delta_1}{8^{-1}b}+ \delta_2})\sqrt{t- \frac{7}{8} b}\Big\}.  \label{dist between pi w and z} 
\end{align}

\begin{figure}[H]
\begin{center}
\includegraphics{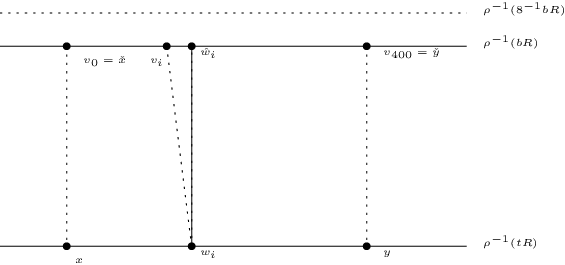}
\caption{The choice of $v_i, w_i, \hat{w}_i$}
\label{figure: thm diam of level sets-vi}
\end{center}
\end{figure}

Therefore from (\ref{assumed Gou-Gu}) (with respect to $8^{-1}b$ instead of $b$) and (\ref{dist between pi w and z}),  we obtain
\begin{align}
|w_iw_{i- 1}|&\leq |\hat{w}_i\hat{w}_{i- 1}|+ \frac{t- b}{8^{-1}b}\cdot \delta_1 bR \nonumber \\
&\leq |v_iv_{i- 1}|+ 40R\Big\{8^{-1}b+\delta_1+ (b+ \sqrt{\frac{\delta_1}{8^{-1}b}+ \delta_2 })\sqrt{t- \frac{7}{8} b}\Big\}.  \nonumber 
\end{align}
Taking the sum with respect to $i$,  using (\ref{sum of zi zi-1}) we obtain
\begin{align}
|w_0w_{400}|&\leq 10^6 R\Big\{b+ \sqrt{\frac{\delta_1}{b}+ \delta_2 }\Big\}+ \sum_{i= 1}^{400} |v_iv_{i- 1}| \nonumber \\
&\leq 10^6 R\Big\{b+ \sqrt{\frac{\delta_1}{b}+ \delta_2 }\Big\}+ 40R(b+ \sqrt{\frac{\delta_1}{b}+ \delta_2})\sqrt{t- 7b} \nonumber \\
&\leq 10^7 R\Big\{b+ \sqrt{\frac{\delta_1}{b}+ \delta_2}\Big\}.  \nonumber 
\end{align}

\textbf{Step (2)}. If $x, y\in \rho^{-1}(sR)$ where $s\in [0, 8b]$, then from $\rho(p)= R$ and the continuity of the function $\rho$, there are 
\begin{align}
\tilde{x}\in \overline{p, x}\cap \rho^{-1}(8bR), \quad \quad \quad \tilde{y}\in \overline{p, y}\cap \rho^{-1}(8bR). \nonumber 
\end{align}

Note 
\begin{align}
|p\tilde{x}|\geq |\rho(p)- \rho(\tilde{x})|= (1- 8b)R. \nonumber 
\end{align}
Then we get 
\begin{align}
|\tilde{x}x|= |px|- |p\tilde{x}|\leq 8bR. \nonumber 
\end{align}
Similar argument yields $|\tilde{y}y|\leq 8bR$. 

Now from Proposition \ref{prop outer level set diam upper bound} and the above, we have
\begin{align}
|xy|&\leq |x\tilde{x}|+ |y\tilde{y}|+ |\tilde{x}\tilde{y}|\leq 16bR+ 80R\Big\{b+ \sqrt{\frac{\delta_1}{b}+ \delta_2}\Big\} \leq 10^2R\Big\{b+ \sqrt{\frac{\delta_1}{b}+ \delta_2}\Big\}, \nonumber 
\end{align}
which yields
\begin{align}
\sup_{x, y\in \rho^{-1}(sR)}|xy|\leq 10^2R\Big\{b+ \sqrt{\frac{\delta_1}{b}+ \delta_2}\Big\}. \nonumber 
\end{align}

For $s\in [2- b, 2]$, the similar argument yields the same upper bound as above. Combining all the above, the conclusion follows.
}
\qed

\section{The width estimate in form of spectral gap}\label{sec main thm}

\begin{theorem}\label{thm width of MCB}
{There is $C(n)> 0$,  such that for any $n$-dim geodesic ball $B_p(R)$ with $Rc\geq 0$ and mean convex boundary,  we have 
\begin{align}
\lambda_1(B_p(R))\cdot R^2- \frac{\pi^2}{16}\geq C(n)\Big\{\frac{\mathcal{W}(B_p(R))}{R}\Big\}^{\frac{4(2n+ 1)(8n+ 15)}{n+ 2}}.   \nonumber 
\end{align}
}
\end{theorem}

\pf
{From Proposition \ref{prop Gou-Gu with spectral gap}, we know that (\ref{assumed Gou-Gu}) holds for 
\begin{align}
\delta_1= C(n)\Big\{b^{-1+ \frac{1}{2n+ 4}}(\epsilon R^2)^{\frac{1}{8(2n+ 1)}}\Big\} . \nonumber 
\end{align}

By Proposition \ref{prop grad curv of rho is extendable}, we know that (\ref{geodesic line extendable}) holds for $\displaystyle \delta_2= C(n)(\frac{\epsilon R^2}{b})^{\frac{1}{n- 1}}$. 

Now by Theorem \ref{thm diam of the level sets}, we get
\begin{align}
&\mathcal{W}(B_p(R))\leq 10^7 R\Big\{b+ \sqrt{\frac{\delta_1 }{b}+ \delta_2}\Big\} \nonumber \\
&\leq C(n)R\cdot \Big(b+ b^{-1+ \frac{1}{4n+ 8}}(\epsilon R^2)^{\frac{1}{16(2n+ 1)}}+ \big(\frac{\epsilon R^2}{b}\big)^{\frac{1}{2(n- 1)}} \Big) \nonumber \\
&\leq C(n)R\cdot \Big(b+ b^{-1+ \frac{1}{4n+ 8}}(\epsilon R^2)^{\frac{1}{16(2n+ 1)}} \Big). \label{the bound of W by b and spectral gap}
\end{align}

Now we choose $b$ such that $\displaystyle b= b^{-1+ \frac{1}{4n+ 8}}(\epsilon R^2)^{\frac{1}{16(2n+ 1)}}$,  which implies
\begin{align}
b= (\epsilon R^2)^{\frac{n+ 2}{4(2n+ 1)(8n+ 15)}}. \nonumber 
\end{align}

Plugging the above choice of $b$ into (\ref{the bound of W by b and spectral gap}), we have
\begin{align}
\mathcal{W}(B_p(R))\leq C(n)R\cdot (\epsilon R^2)^{\frac{n+ 2}{4(2n+ 1)(8n+ 15)}}= C(n)R\cdot (\lambda R^2- \frac{\pi^2}{16})^{\frac{n+ 2}{4(2n+ 1)(8n+ 15)}}. \nonumber 
\end{align}
The conclusion follows.
}
\qed

\section*{Acknowledgments}
We thank Zuoqin Wang for discussion and helpful suggestion. And we are grateful to an anonymous referee for his detailed suggestion and comments on the earlier version of this paper, which greatly improves this paper. 

\appendix
\section{A generalization of the first variation formula}
\begin{center}
    By Zichang Liu\footnote{Department of Mathematical Sciences, Tsinghua University, Beijing, P. R. China. E-mail address: liu-zc19@mails.tsinghua.edu.cn}
\end{center}

\begin{theorem}\label{thm general of 1st vari form}

Assume $M$ is an arbitrary connected Riemannian manifold, and $\tau_1(s),\ \tau_2(s):(0,T)\rightarrow M$ are (possibly constant) geodesics on $M$. Let $|\cdot,\cdot|$ be the distance function on $M$. Assume\\
(i) $\rho(s)=|\tau_1(s),\tau_2(s)|$ is differentiable at $s=s_0$;\\
(ii) $\tau_1(s_0)\neq \tau_2(s_0)$, and there exists a minimal (unit speed) geodesic $\phi$ joining $\tau_1(s_0)$ to $\tau_2(s_0)$ (with domain $[0,\rho(s_0)]$).\\
Then
$$\frac{d\rho(s)}{ds}|_{s=s_0}=\langle \phi'(\rho(s_0)),\tau_2'(s_0)\rangle-\langle \phi'(0),\tau_1'(s_0)\rangle.$$
\end{theorem}
\proof
Let $V_1(t)$, $V_2(t)$ be parallel vector fields along $\phi$ with 
$$V_1(0)=\tau_1'(s_0),\ \ V_2(\rho(s_0))=\tau_2'(s_0),$$
and construct a family of smooth curves:
$$\Phi(t,s)=\exp_{\phi(t)}[(s-s_0)((1-\frac{t}{\rho(s_0)})V_1(t)+\frac{t}{\rho(s_0)}V_2(t))],\ \ t\in[0,\rho(s_0)],\ \  s\in(s_0-\epsilon,s_0+\epsilon).$$
By the first variation formula,
$$\frac{dl(\Phi(\cdot,s))}{ds}|_{s=s_0}=\langle \phi'(\rho(s_0)),\tau_2'(s_0)\rangle-\langle \phi'(0),\tau_1'(s_0)\rangle.$$
Since $\rho(s)$ is differentiable at $s_0$, noticing the fact that
$$\Phi(t,s_0)=\phi(t),\ \ l(\Phi(\cdot,s_0))=\rho(s_0),\ \ l(\Phi(\cdot,s))\geq \rho(s),\ \forall s,$$
we get the result.
\qed

When one of the two geodesics $\tau_1,\ \tau_2$ degenerates to a point, we get the following conclusion:
\begin{cor}
Assume $M$ is an arbitrary connected Riemannian manifold, $x\in M$, and $\tau(s):(0,T)\rightarrow M$ is a (unit speed) geodesic on $M$. Assume\\
(i) $\rho(s)=|\tau(s),x|$ is differentiable at $s=s_0$;\\
(ii) there exists a minimal (unit speed) geodesic $\phi$ joining $x$ to $\tau(s_0)$ (with domain $[0,\rho(s_0)]$).\\
Then
$$\frac{d\rho(s)}{ds}|_{s=s_0}=\cos\angle(\tau'(s_0),\phi'(\rho(s_0))).$$
\end{cor}


\begin{bibdiv}
\begin{biblist}


\bib{Aubry}{article}{
    AUTHOR = {Aubry, Erwann},
     TITLE = {Pincement sur le spectre et le volume en courbure de {R}icci
              positive},
   JOURNAL = {Ann. Sci. \'{E}cole Norm. Sup. (4)},
  FJOURNAL = {Annales Scientifiques de l'\'{E}cole Normale Sup\'{e}rieure. Quatri\`eme
              S\'{e}rie},
    VOLUME = {38},
      YEAR = {2005},
    NUMBER = {3},
     PAGES = {387--405},
      ISSN = {0012-9593},
   MRCLASS = {53C21 (53C20)},
  MRNUMBER = {2166339},
MRREVIEWER = {Nader Yeganefar},
       DOI = {10.1016/j.ansens.2005.01.002},
       URL = {https://doi.org/10.1016/j.ansens.2005.01.002},
}

\bib{BDV}{article}{
    AUTHOR = {Brasco, Lorenzo},
    author= {De Philippis, Guido },
    author= {Velichkov, Bozhidar},
     TITLE = {Faber-{K}rahn inequalities in sharp quantitative form},
   JOURNAL = {Duke Math. J.},
  FJOURNAL = {Duke Mathematical Journal},
    VOLUME = {164},
      YEAR = {2015},
    NUMBER = {9},
     PAGES = {1777--1831},
      ISSN = {0012-7094},
   MRCLASS = {49R05 (47A75 49Q20)},
  MRNUMBER = {3357184},
MRREVIEWER = {Antoine Lemenant},
       DOI = {10.1215/00127094-3120167},
       URL = {https://doi.org/10.1215/00127094-3120167},
}

\bib{Buser}{article}{
    AUTHOR = {Buser, Peter},
     TITLE = {A note on the isoperimetric constant},
   JOURNAL = {Ann. Sci. \'Ecole Norm. Sup. (4)},
  FJOURNAL = {Annales Scientifiques de l'\'Ecole Normale Sup\'erieure.
              Quatri\`eme S\'erie},
    VOLUME = {15},
      YEAR = {1982},
    NUMBER = {2},
     PAGES = {213--230},
      ISSN = {0012-9593},
     CODEN = {ASENAH},
   MRCLASS = {58G25 (52A40 53C20)},
  MRNUMBER = {683635 (84e:58076)},
MRREVIEWER = {Scott Wolpert},
       URL = {http://www.numdam.org/item?id=ASENS_1982_4_15_2_213_0},
}

\bib{CMM}{article}{
    AUTHOR = {Cavalletti, F.},
    author= {Maggi, F.},
    author= {Mondino, A.},
     TITLE = {Quantitative isoperimetry \`a la {L}evy-{G}romov},
   JOURNAL = {Comm. Pure Appl. Math.},
  FJOURNAL = {Communications on Pure and Applied Mathematics},
    VOLUME = {72},
      YEAR = {2019},
    NUMBER = {8},
     PAGES = {1631--1677},
      ISSN = {0010-3640},
   MRCLASS = {53C42 (49Q15 53C23)},
  MRNUMBER = {3974951},
MRREVIEWER = {Luis Guijarro},
       DOI = {10.1002/cpa.21808},
       URL = {https://doi.org/10.1002/cpa.21808},
}
	
\bib{CMS}{article}{
    AUTHOR = {Cavalletti, Fabio},
    author= {Mondino, Andrea},
    author= {Semola, Daniele},
     TITLE = {Quantitative {O}bata's theorem},
   JOURNAL = {Anal. PDE},
  FJOURNAL = {Analysis \& PDE},
    VOLUME = {16},
      YEAR = {2023},
    NUMBER = {6},
     PAGES = {1389--1431},
      ISSN = {2157-5045},
   MRCLASS = {58J50 (53C23)},
  MRNUMBER = {4632377},
       DOI = {10.2140/apde.2023.16.1389},
       URL = {https://doi.org/10.2140/apde.2023.16.1389},
}

\bib{Chavel}{book}{
    AUTHOR = {Chavel, Isaac},
     TITLE = {Eigenvalues in {R}iemannian geometry},
    SERIES = {Pure and Applied Mathematics},
    VOLUME = {115},
      NOTE = {Including a chapter by Burton Randol,
              With an appendix by Jozef Dodziuk},
 PUBLISHER = {Academic Press, Inc., Orlando, FL},
      YEAR = {1984},
     PAGES = {xiv+362},
      ISBN = {0-12-170640-0},
   MRCLASS = {58G25 (35P99 53C20)},
  MRNUMBER = {768584},
MRREVIEWER = {G\'{e}rard Besson},
}

\bib{CC-Ann}{article}{
    AUTHOR = {Cheeger, Jeff},
    author= {Colding, Tobias H.},
     TITLE = {Lower bounds on {R}icci curvature and the almost rigidity of
              warped products},
   JOURNAL = {Ann. of Math. (2)},
  FJOURNAL = {Annals of Mathematics. Second Series},
    VOLUME = {144},
      YEAR = {1996},
    NUMBER = {1},
     PAGES = {189--237},
      ISSN = {0003-486X},
     CODEN = {ANMAAH},
   MRCLASS = {53C21 (53C20 53C23)},
  MRNUMBER = {1405949 (97h:53038)},
MRREVIEWER = {Joseph E. Borzellino},
       DOI = {10.2307/2118589},
       URL = {http://dx.doi.org/10.2307/2118589},
}

\bib{Colding-sphere-1}{article}{
    AUTHOR = {Colding, Tobias H.},
     TITLE = {Shape of manifolds with positive {R}icci curvature},
   JOURNAL = {Invent. Math.},
  FJOURNAL = {Inventiones Mathematicae},
    VOLUME = {124},
      YEAR = {1996},
    NUMBER = {1-3},
     PAGES = {175--191},
      ISSN = {0020-9910},
   MRCLASS = {53C23 (53C21)},
  MRNUMBER = {1369414},
MRREVIEWER = {Man Chun Leung},
       DOI = {10.1007/s002220050049},
       URL = {https://doi.org/10.1007/s002220050049},
}

\bib{Colding-sphere-2}{article}{
    AUTHOR = {Colding, Tobias H.},
     TITLE = {Large manifolds with positive {R}icci curvature},
   JOURNAL = {Invent. Math.},
  FJOURNAL = {Inventiones Mathematicae},
    VOLUME = {124},
      YEAR = {1996},
    NUMBER = {1-3},
     PAGES = {193--214},
      ISSN = {0020-9910},
   MRCLASS = {53C23 (53C21)},
  MRNUMBER = {1369415},
MRREVIEWER = {Man Chun Leung},
       DOI = {10.1007/s002220050050},
       URL = {https://doi.org/10.1007/s002220050050},
}

\bib{CN}{article}{
    AUTHOR = {Colding, Tobias Holck},
    author= {Naber, Aaron},
     TITLE = {Sharp {H}\"older continuity of tangent cones for spaces with a lower {R}icci curvature bound and applications},
   JOURNAL = {Ann. of Math. (2)},
  FJOURNAL = {Annals of Mathematics. Second Series},
    VOLUME = {176},
      YEAR = {2012},
    NUMBER = {2},
     PAGES = {1173--1229},
      ISSN = {0003-486X},
     CODEN = {ANMAAH},
   MRCLASS = {53C21 (53C20)},
  MRNUMBER = {2950772},
MRREVIEWER = {Yu Ding},
       DOI = {10.4007/annals.2012.176.2.10},
       URL = {http://dx.doi.org/10.4007/annals.2012.176.2.10},
}

\bib{FMP}{article}{
    AUTHOR = {Fusco, N.},
    author= {Maggi, F.},
    author= {Pratelli, A.},
     TITLE = {The sharp quantitative isoperimetric inequality},
   JOURNAL = {Ann. of Math. (2)},
  FJOURNAL = {Annals of Mathematics. Second Series},
    VOLUME = {168},
      YEAR = {2008},
    NUMBER = {3},
     PAGES = {941--980},
      ISSN = {0003-486X},
   MRCLASS = {52A40 (28A99)},
  MRNUMBER = {2456887},
MRREVIEWER = {Sasha Sodin},
       DOI = {10.4007/annals.2008.168.941},
       URL = {https://doi.org/10.4007/annals.2008.168.941},
}

\bib{HW}{article}{
    AUTHOR = {Hang, Fengbo},
    author= {Wang, Xiaodong},
     TITLE = {A remark on {Z}hong-{Y}ang's eigenvalue estimate},
   JOURNAL = {Int. Math. Res. Not. IMRN},
  FJOURNAL = {International Mathematics Research Notices. IMRN},
      YEAR = {2007},
    NUMBER = {18},
     PAGES = {Art. ID rnm064, 9},
      ISSN = {1073-7928},
   MRCLASS = {53C21 (58J50)},
  MRNUMBER = {2358887},
MRREVIEWER = {Fr\'{e}d\'{e}ric Robert},
       DOI = {10.1093/imrn/rnm064},
       URL = {https://doi.org/10.1093/imrn/rnm064},
}

\bib{Hersch}{article}{
    AUTHOR = {Hersch, Joseph},
     TITLE = {Sur la fr\'{e}quence fondamentale d'une membrane vibrante:
              \'{e}valuations par d\'{e}faut et principe de maximum},
   JOURNAL = {Z. Angew. Math. Phys.},
  FJOURNAL = {Zeitschrift f\"{u}r Angewandte Mathematik und Physik. ZAMP.
              Journal of Applied Mathematics and Physics. Journal de
              Math\'{e}matiques et de Physique Appliqu\'{e}es},
    VOLUME = {11},
      YEAR = {1960},
     PAGES = {387--413},
      ISSN = {0044-2275},
   MRCLASS = {35.80},
  MRNUMBER = {125319},
MRREVIEWER = {H. F. Weinberger},
       DOI = {10.1007/BF01604498},
       URL = {https://doi-org-s.qh.yitlink.com:8444/10.1007/BF01604498},
}

\bib{Kasue}{article}{
    AUTHOR = {Kasue, Atsushi},
     TITLE = {On a lower bound for the first eigenvalue of the {L}aplace
              operator on a {R}iemannian manifold},
   JOURNAL = {Ann. Sci. \'{E}cole Norm. Sup. (4)},
  FJOURNAL = {Annales Scientifiques de l'\'{E}cole Normale Sup\'{e}rieure. Quatri\`eme
              S\'{e}rie},
    VOLUME = {17},
      YEAR = {1984},
    NUMBER = {1},
     PAGES = {31--44},
      ISSN = {0012-9593},
   MRCLASS = {58G25 (35P15)},
  MRNUMBER = {744066},
MRREVIEWER = {G. Tsagas},
       URL = {https://www-numdam-org-s.qh.yitlink.com:8444/item?id=ASENS_1984_4_17_1_31_0},
}

\bib{Klartag}{article}{
    AUTHOR = {Klartag, Bo'az},
     TITLE = {Needle decompositions in {R}iemannian geometry},
   JOURNAL = {Mem. Amer. Math. Soc.},
  FJOURNAL = {Memoirs of the American Mathematical Society},
    VOLUME = {249},
      YEAR = {2017},
    NUMBER = {1180},
     PAGES = {v+77},
      ISSN = {0065-9266},
      ISBN = {978-1-4704-2542-5; 978-1-4704-4127-2},
   MRCLASS = {53C21 (52A20 52A40)},
  MRNUMBER = {3709716},
MRREVIEWER = {Vasyl Gorkavyy},
       DOI = {10.1090/memo/1180},
       URL = {https://doi.org/10.1090/memo/1180},
}

\bib{LiBook}{book}{
   author={Li, Peter},
   title={Geometric analysis},
   place={Cambridge Studies in Advanced Mathematics, 134. Cambridge University Press, Cambridge, x+406 pp},
    date={2012}, 
   }

\bib{LY}{incollection}{
    AUTHOR = {Li, Peter},
    author = {Yau, Shing Tung},
     TITLE = {Estimates of eigenvalues of a compact {R}iemannian manifold},
 BOOKTITLE = {Geometry of the {L}aplace operator ({P}roc. {S}ympos. {P}ure
              {M}ath., {U}niv. {H}awaii, {H}onolulu, {H}awaii, 1979)},
    SERIES = {Proc. Sympos. Pure Math., XXXVI},
     PAGES = {205--239},
 PUBLISHER = {Amer. Math. Soc., Providence, R.I.},
      YEAR = {1980},
   MRCLASS = {58G25 (53C20)},
  MRNUMBER = {573435},
MRREVIEWER = {P. G\"{u}nther},
}

\bib{Liu}{article}{
    AUTHOR = {Liu, Zichang},
     TITLE = {A generalization of the first variation formula},
   JOURNAL = {Preprint},
}

\bib{MH}{article}{
    AUTHOR = {M\'{e}ndez-Hern\'{a}ndez, Pedro J.},
     TITLE = {Brascamp-{L}ieb-{L}uttinger inequalities for convex domains of
              finite inradius},
   JOURNAL = {Duke Math. J.},
  FJOURNAL = {Duke Mathematical Journal},
    VOLUME = {113},
      YEAR = {2002},
    NUMBER = {1},
     PAGES = {93--131},
      ISSN = {0012-7094},
   MRCLASS = {31B35},
  MRNUMBER = {1905393},
MRREVIEWER = {Catherine Bandle},
       DOI = {10.1215/S0012-7094-02-11313-1},
       URL = {https://doi.org/10.1215/S0012-7094-02-11313-1},
}

\bib{PW}{article}{
    AUTHOR = {Payne, L. E.},
    author = {Weinberger, H. F.},
     TITLE = {An optimal {P}oincar\'{e} inequality for convex domains},
   JOURNAL = {Arch. Rational Mech. Anal.},
  FJOURNAL = {Archive for Rational Mechanics and Analysis},
    VOLUME = {5},
      YEAR = {1960},
     PAGES = {286--292 (1960)},
      ISSN = {0003-9527},
   MRCLASS = {35.00},
  MRNUMBER = {117419},
MRREVIEWER = {I. Stakgold},
       DOI = {10.1007/BF00252910},
       URL = {https://doi.org/10.1007/BF00252910},
}

\bib{Petersen}{article}{
    AUTHOR = {Petersen, Peter},
     TITLE = {On eigenvalue pinching in positive {R}icci curvature},
   JOURNAL = {Invent. Math.},
  FJOURNAL = {Inventiones Mathematicae},
    VOLUME = {138},
      YEAR = {1999},
    NUMBER = {1},
     PAGES = {1--21},
      ISSN = {0020-9910},
   MRCLASS = {53C20 (53C21 58J50)},
  MRNUMBER = {1714334},
MRREVIEWER = {Joseph E. Borzellino},
       DOI = {10.1007/s002220050339},
       URL = {https://doi.org/10.1007/s002220050339},
}

\bib{Protter}{article}{
    AUTHOR = {Protter, M. H.},
     TITLE = {A lower bound for the fundamental frequency of a convex
              region},
   JOURNAL = {Proc. Amer. Math. Soc.},
  FJOURNAL = {Proceedings of the American Mathematical Society},
    VOLUME = {81},
      YEAR = {1981},
    NUMBER = {1},
     PAGES = {65--70},
      ISSN = {0002-9939},
   MRCLASS = {35P15},
  MRNUMBER = {589137},
MRREVIEWER = {J. Hersch},
       DOI = {10.2307/2043987},
       URL = {https://doi-org-s.qh.yitlink.com:8444/10.2307/2043987},
}

\bib{WXZ}{article}{
    AUTHOR = {Wang, Haibin},
    author = {Xu, Guoyi},
    author = {Zhou, Jie},
     TITLE = {The sharp lower bound of the first {D}irichlet eigenvalue for
              geodesic balls},
   JOURNAL = {Math. Z.},
  FJOURNAL = {Mathematische Zeitschrift},
    VOLUME = {300},
      YEAR = {2022},
    NUMBER = {2},
     PAGES = {2063--2068},
      ISSN = {0025-5874},
   MRCLASS = {58J50},
  MRNUMBER = {4363806},
       DOI = {10.1007/s00209-021-02859-8},
       URL = {https://doi.org/10.1007/s00209-021-02859-8},
}
  
\bib{Xu}{article}{
    AUTHOR = {Xu, Guoyi},
     TITLE = {Local estimate of fundamental groups},
   JOURNAL = {Adv. Math.},
  FJOURNAL = {Advances in Mathematics},
    VOLUME = {352},
      YEAR = {2019},
     PAGES = {158--230},
      ISSN = {0001-8708},
   MRCLASS = {53C21 (57M05)},
  MRNUMBER = {3959654},
MRREVIEWER = {Christine M. Escher},
       DOI = {10.1016/j.aim.2019.06.006},
       URL = {https://doi-org-s.qh.yitlink.com:8444/10.1016/j.aim.2019.06.006},
} 

\bib{XZ}{article}{
    AUTHOR = {Xu, Guoyi},
    author= {Zhou, Jie},
     TITLE = {The construction of {$\epsilon$}-splitting map},
   JOURNAL = {Calc. Var. Partial Differential Equations},
  FJOURNAL = {Calculus of Variations and Partial Differential Equations},
    VOLUME = {62},
      YEAR = {2023},
    NUMBER = {3},
     PAGES = {Paper No. 75},
      ISSN = {0944-2669},
   MRCLASS = {53 (35K15)},
  MRNUMBER = {4531769},
       DOI = {10.1007/s00526-022-02418-x},
       URL = {https://doi-org-s.qh.yitlink.com:8444/10.1007/s00526-022-02418-x},
} 

\bib{WZ}{article}{
    AUTHOR = {Wang, Bing},
    author= {Zhao, Xinrui},
     TITLE = {Canonical diffeomorphisms of manifolds near spheres},
   JOURNAL = {J. Geom. Anal.},
  FJOURNAL = {Journal of Geometric Analysis},
    VOLUME = {33},
      YEAR = {2023},
    NUMBER = {9},
     PAGES = {Paper No. 304, 31},
      ISSN = {1050-6926,1559-002X},
   MRCLASS = {53C20},
  MRNUMBER = {4615511},
       DOI = {10.1007/s12220-023-01375-x},
       URL = {https://doi-org-s.qh.yitlink.com:8444/10.1007/s12220-023-01375-x},
}

\bib{Yang}{article}{
    AUTHOR = {Yang, DaGang},
     TITLE = {Lower bound estimates of the first eigenvalue for compact
              manifolds with positive {R}icci curvature},
   JOURNAL = {Pacific J. Math.},
  FJOURNAL = {Pacific Journal of Mathematics},
    VOLUME = {190},
      YEAR = {1999},
    NUMBER = {2},
     PAGES = {383--398},
      ISSN = {0030-8730},
   MRCLASS = {53C21 (53C20 58J50)},
  MRNUMBER = {1722898},
MRREVIEWER = {William P. Minicozzi, II},
       DOI = {10.2140/pjm.1999.190.383},
       URL = {https://doi.org/10.2140/pjm.1999.190.383},
}

\bib{ZY}{article}{
    AUTHOR = {Zhong, Jia Qing},
    author = {Yang, Hong Cang},
     TITLE = {On the estimate of the first eigenvalue of a compact
              {R}iemannian manifold},
   JOURNAL = {Sci. Sinica Ser. A},
  FJOURNAL = {Scientia Sinica. Series A. Mathematical, Physical,
              Astronomical \& Technical Sciences},
    VOLUME = {27},
      YEAR = {1984},
    NUMBER = {12},
     PAGES = {1265--1273},
      ISSN = {0253-5831},
   MRCLASS = {58G25},
  MRNUMBER = {794292},
MRREVIEWER = {Domenico Perrone},
}
	
\end{biblist}
\end{bibdiv}
\end{document}